\documentclass[11pt,a4paper]{article}
\usepackage{amsmath}
\usepackage{amsfonts}
\usepackage{amssymb}
\usepackage{amsthm}
\usepackage{authblk}
\usepackage{mathrsfs}
\usepackage{dsfont}
\usepackage{paralist}
\usepackage{natbib}
\usepackage{bm}
\usepackage{bbm}
\usepackage{color}
\usepackage{graphicx}
\usepackage{tabularx}
\usepackage[left=3.25cm,right=3.25cm,top=2cm,bottom=2cm,includeheadfoot]{geometry}
\usepackage{comment}
\usepackage{tikz}
\usetikzlibrary{decorations.pathreplacing}
\usetikzlibrary{decorations.pathmorphing}
\usepackage{resizegather}
\DeclareFontFamily{OT1}{pzc}{}
\DeclareFontShape{OT1}{pzc}{m}{it}{<-> s * [1.10] pzcmi7t}{}
\DeclareMathAlphabet{\mathpzc}{OT1}{pzc}{m}{it}
\usepackage[colorlinks=true,linkcolor=blue,citecolor=blue,pdfborder={0 0 0}]{hyperref}

\newtheoremstyle{normal}
{2ex}               
{3ex}               
{}                  
{}                  
{\bfseries} 
{}                  
{2pt}   
{\thmname{#1}\thmnumber{ #2.} \thmnote{(#3)}}

\newtheoremstyle{italic}
{2ex}
{3ex}
{\itshape}
{}
{\bfseries} 
{}
{2pt}
{\thmname{#1}\thmnumber{ #2.} \thmnote{(#3)}}

\theoremstyle{normal}
\newtheorem{definition}{Definition}[section]
\newtheorem{remark}[definition]{Remark}

\newtheorem{condition}[definition]{Condition}

\theoremstyle{italic}
\newtheorem{theorem}[definition]{Theorem}
\newtheorem{lemma}[definition]{Lemma}

\renewcommand{\P}{\mathbb{P}}

\newcommand{\E}{\mathbb{E}}

\DeclareMathAccent{\verywidehat}{\mathord}{largesymbols}{'144}

\newcommand{\R}{\mathbb{R}}

\newcommand{\al}{\alpha}
\newcommand{\be}{\beta}
\newcommand{\la}{\lambda}
\newcommand{\La}{\Lambda}
\newcommand{\ga}{\gamma}

\newcommand{\Ga}{\Gamma}
\newcommand{\si}{\sigma}

\newcommand{\vpi}{\varpi}

\newcommand{\de}{\delta}

\newcommand{\De}{\Delta}

\newcommand{\Om}{\Omega}
\newcommand{\ze}{\zeta}

\newcommand{\nN}{\mathcal N}

\newcommand{\lL}{\mathcal L}

\newcommand{\tols}{~\stackrel{\lL-(s)}{\longrightarrow}~}
\newcommand{\tol}{\stackrel{\lL}{\longrightarrow}}

\newcommand{\pn}{\stackrel{\P}{\longrightarrow}}

\DeclareMathOperator{\Var}{Var}

\allowdisplaybreaks[1]

\begin{document}

\title{Jump regression revisited}
\author{Mathias Vetter\thanks{corresponding author; postal address: Christian-Albrechts-Universit\"at zu Kiel, Mathematisches Seminar, Heinrich-Hecht-Platz\ 6, 24118 Kiel, Germany; e-mail address: vetter@math.uni-kiel.de.} }
\author{Fan Yu}
\affil{Christian-Albrechts-Universit\"at zu Kiel}


\maketitle

\begin{abstract}
Regression analysis for stochastic processes has been an important topic over the last decades. Typically, regressing a dependent process $Y$ on an explanatory $Z$ results in an integral relationship of the form
$dY_t = \be_t dZ_t + dX_t$ with some residual process $X$, and often with constant or piecewise constant $\be$. In the case of jump processes such a relationship essentially boils down to $\De Y_t = \be \De Z_t + \De X_t$ on the level of the jumps. Here the residual process $X$ is such that it never jumps together with $Z$, and so any jump in $Y$ is either exactly proportional to $Z$ or purely idiosyncratic. In this paper we discuss a related setup which appears more realistic and is closer to classical regression, namely $\De Y_t = (\be + \eta_t) \De Z_t$ for i.i.d.\ $(\eta_t)_{t \ge 0}$. We propose an estimator for $\be$ based on high-frequency observations of $(Y,Z)$ and discuss consistency and associated central limit theorems in two different asymptotic regimes. 
\end{abstract}

\medskip

\textit{Keywords and Phrases:} Central limit theorem, compound Poisson process, jump regression, regression analysis, stable process  

\smallskip

\textit{AMS Subject Classification:} 62F10, 62M09, 62M10


\section{Introduction} \label{sec:int}
\def\theequation{1.\arabic{equation}}
\setcounter{equation}{0}

Understanding dependence between variables is not only a classical problem in finite-dimensional statistics but has gained at lot of attention in recent years in the context of stochastic processes as well. Here the setup usually is as follows: The dependent process $Y$ satisfies 
\begin{align} \label{defY}
Y_t = Y_0 + \int_0^t \be_s dZ_s + X_t
\end{align} 
with a covariate process $Z$ and a residual process $X$. Based on observations of $(Y,Z)$ the goal is either to estimate the process $\be$ (which in some models only is a constant) or to test certain hypotheses about it. Typical applications are in the context of finance where \(Z\) often serves as a risk factor process, e.g., a market index or a macroeconomic indicator.

Early examples in this context are \cite{barshe2004} and \cite{mykzha2006}, among others, who deal with continuous semimartingales only and generalize standard regression techniques to the setting of stochastic processes observed at high-frequency. Ultimately, if the process $(Y,Z)$ was observed continuously, the standard estimators for $\be$ would be related to certain covariation processes, and hence the statistical techniques heavily rely on its discrete time analogues such as realized variance and realized covariance. 

Since then the attention has shifted slightly towards factor models which are not only based on continuous processes but involve jumps in $Y$ and $Z$ as well. This includes statistical situations where the jumps essentially are just a nuisance parameter (see e.g.\ \citealp{lietal2016}), but the dependence between jumps has gained attention as well. \cite{todbol2010} discuss a regression model in which the jumps and the diffusive part may have different (but constant) betas while \cite{aitetal2020} discuss the estimation of a time-varying beta in a similar situation. See also \cite{boletal2016} and \cite{aitetal2020a} for the estimation of risk premia in continuous-time factor models involving jumps. 

The starting point for our paper is a series of papers by Li and co-authors where the focus is on regression between jumps only. Specifically, if the relationship between $Y$ and $Z$ as in (\ref{defY}) only holds on the level of the jumps, then this translates into
\begin{align} \label{defju}
\De Y_t = \be_t \De Z_t + \De X_t,
\end{align}
where the residual process $X$ is such that it never jumps together with $Z$. Hence, any jump in $Y$ is either exactly proportional to a jump in $Z$ or purely idiosyncratic (i.e.\, coming from $X$). This is the essence of the model in \cite{lietal2017a}, \cite{lietal2017b} and \cite{lietal2017c}, in which the estimation of $\be$ and a test of temporal constancy is discussed when $Y$ contains additional diffusive components. It is in particular this additional Brownian part which allows for an asymptotic theory in their works. 

There is a downside with the model (\ref{defju}) and, in turn, (\ref{defY}), however. At least on the level of the jumps, the assumption of an exact proportionality between $\De Y_t$ and $\De Z_t$ is unrealistic from a practical point of view, in particular in the typical setup of a constant (or piecewise constant) $\be$. This insight was the motivation for the statistical analysis in this paper. 

As discussed below, the key change in our statistical model is that we allow for a connection between the jumps of $Y$ and the jumps of $Z$ which is of the form 
\begin{align*} 
  \Delta Y_t = (\beta+\eta_t) \Delta Z_t = \beta \De Z_t + \eta_t \Delta Z_t
\end{align*}
for i.i.d.\ error variables $\eta_t$ with mean zero. Here, for simplicity with a constant $\be$, the jumps of both processes are now only proportional \emph{on average}, resembling much closer the classical situation in linear regression. Our paper is then theoretical by nature. We are first and foremost interested in the question how an estimation of $\be$ could work in this new setting. Here we discuss different observational situations for which in turn different assumptions on the regressor $Z$ are needed.   

The paper is organised as follows: In Section \ref{sec:set} the basic setup in this paper is presented and discussed, while all theoretical results are given in Section \ref{sec:res}. Here we deal both with a finite and an infinite time horizon situation. A short simulation study is given in Section \ref{sec:sim}. All proofs are gathered in Section \ref{sec:proof}.

\section{Setting} \label{sec:set}
\def\theequation{2.\arabic{equation}}
\setcounter{equation}{0}
We assume a probability space \( (\Omega,\mathscr{F},\P) \) with a filtration \((\mathscr{F}_t)_{t \ge 0}\) that satisfies the usual conditions. On this space we consider a bivariate Itô semimartingale \((Y,Z)\) involving jumps, defined via
\begin{align}  \label{defYZ}
  &Y_t=Y_0+ \int_0^t a_s ds+\int_0^t \sigma_s dW_s+\sum_{s\le t} \Delta Y_s, \nonumber \\ 
  &Z_t = Z_0+ \int_0^t a'_s ds+\int_0^t \sigma'_s dW'_s+\sum_{s\le t} \Delta Z_s,\quad t \ge 0.
\end{align}
Here, \( W \) and \( W' \) represent two independent standard Brownian motions, the initial values \( (Y_0, Z_0) \) are \( \mathscr{F}_0 \)-measurable, and the processes $a, a', \si, \si'$ are all assumed to be locally bounded and c\`adl\`ag. 

As we are interested in jump regressions we have to be more specific about the nature of the jump parts of the dependent process $Y$ and the regressor $Z$. For the process $Z$ we assume jumps of finite variation, i.e.\ summable jumps as in (\ref{defYZ}), whereas for the jump process of $Y$ we work with a 
model which resembles linear regression. Precisely, we set
\begin{align} \label{defJ}
  \Delta Y_t = (\beta+\eta_t) \Delta Z_t = \beta \De Z_t + \eta_t \Delta Z_t
\end{align}
where $\be$ is an unknown real parameter and \(\eta = (\eta_t)_{t \ge 0}\) is a sequence of i.i.d.\ bounded random variables with zero mean and finite variance \(\sigma^2_\eta\), adapted to the filtration  \((\mathscr{F}_t)_{t \ge 0}\) and independent of all other processes in the model. 

Throughout this work the goal is to estimate the unknown parameter $\be$ based on high-frequency observations of $(Y,Z)$. Thereby we will work in two different situations regarding the nature of the time horizon and, hence, regarding the regressor $Z$. 

\begin{condition}  \label{condModel}
We assume that the process $(Y,Z)$ is observed at discrete time points $i \De_n$, $0 \le i \le n$, for some sequence $\De_n \longrightarrow 0$, with the number of observations being $n+1$ (and, hence, the number of increments being $n$). 
\begin{itemize}
	\item[(i)] In the \emph{finite time horizon} case we assume that $n = \De_n^{-1}$. Hence, we are dealing with high-frequency observations over the finite interval $[0,1]$. 
	\item[(ii)] In the \emph{infinite time horizon} case we assume that $T_n = n \De_n \longrightarrow \infty$. Thus, we are having high-frequency observations over a potentially infinitely long time interval.
\end{itemize}
\end{condition}

\begin{remark}
\begin{itemize}
	\item[(i)] The jump regression model in \eqref{defJ} is just one of many reasonable definitions, and one could alternatively work with $\Delta Y_t = \beta(1+\eta_t) \Delta Z_t$ as well for which similar results could be derived. What is natural in both cases is to assume that the residual process somehow involves independent noise, but we also assume a scaling property in the sense that the jump size of $Y$ then is essentially proportional to the jump size of $Z$. This rules out a model such as $\Delta Y_t = \beta \Delta Z_t + \eta_t$ which looks like the classical regression case, but means that very small jumps in $Z$ are often transformed to large shocks in $Y$. This might not be very desirable from an application point of view on one hand, but also poses severe theoretical problems as e.g.\ summability of the jumps of $Z$ does then in general not translate to summability of the jumps of $Y$.
	\item[(ii)] While the authors in \cite{lietal2017a} work with finite activity jumps over a finite time interval, such an assumption would not allow for a consistent estimation of $\beta$ in our setup. The crucial difference between the two models is that they assume the jumps of $Y$ to be \emph{exactly proportional} to the jumps of $Z$--an extremely simplifying condition as asymptotically already one observation of a pair of jumps allows to read off the value of $\beta$. In our model infinitely many jumps are necessary as otherwise the influence of $\eta$ cannot be neglected asymptotically. Hence, at least in the finite time horizon case, a model with infinitely many jumps of $Z$ is necessary. 
\end{itemize}
\end{remark}

\section{Main results} \label{sec:res}
\def\theequation{3.\arabic{equation}}
\setcounter{equation}{0}
In this chapter we discuss the asymptotic results for an estimator of $\beta$ in both the finite and the infinite time horizon case. Due to the different nature of the two setups different assumptions on $(Y,Z)$ and the estimator are needed. The basic structure of the estimator is the same, however: Let \(K_n =\Delta_n^\varpi \) for some $0<\varpi <1$. Set also \(P_n = \sum_{i=1}^n \mathbbm{1}_{\{ |\Delta_i^n Z| > K_n \}}\) as the count of increments exceeding the threshold \(K_n\), where we use the shorthand $\De_i^n X = X_{i \De_n} - X_{(i-1)\De_n}$. Then we define 
\begin{align}
  \label{eq:estimator_def}
  \hat{\beta}_n =
  \begin{cases}
    \frac{1}{P_n} \sum_{i=1}^n  \frac{\Delta_i^n Y}{\Delta_i^n Z} \mathbbm{1}_{\left\{ |\Delta_i^n Z| > K_n \right\}} & \text{if } P_n > 0, \\
    0 & \text{if } P_n = 0.
  \end{cases}
\end{align}

\begin{remark}
\begin{itemize}
	\item[(i)] 
It is natural to base the estimation of $\be$ on the ratio between ${\Delta_i^n Y}$ and ${\Delta_i^n Z}$, as the unknown parameter resembles the average ratio of the corresponding jump sizes of $Y$ and $Z$. Working with a threshold of $|\Delta_i^n Z| > K_n$ guarantees that only those incrememts of $Z$ are taken into account which are reasonably large and, hence, most likely dominated by a single jump and not by drift and martingale parts. See e.g.\ \cite{mancini2009} for more background on thresholding techniques in high-frequency statistics. $\hat \be_n$ is then simply defined as the average over such ratios.
\item[(ii)] It is in principle possible to work with a weighted version of $\hat \be_n$ as well. Our goal, however, is to understand how estimation of $\be$ in our model works in principle and not to obtain optimal rates of convergences or even variances. Note also that the standard choice of a threshold is typically of the form $K_n = C \De_n^\varpi$ for some $C > 0$. Such a variant is possible in our setting as well, but might--depending on the statistical setting--change the limiting variance in the central limit theorem.  
\end{itemize}
\end{remark}

\subsection{The finite time horizon} \label{secfin}
In this section we assume high-frequency observations of $(Y,Z)$ over the interval $[0,1]$, i.e., we have data $(Y_{i\De_n}, Z_{i \De_n})$ available with the number of increments $n= \De_n^{-1}$. Before we get to the main result in this section, let us recall a classical definition from probability theory. 

\begin{definition} \label{defstab}
\begin{itemize}
	\item[(i)] A real-valued random variable $X$ is called \emph{stable} if its characteristic function $\Phi(u) = \E[\exp(iuX)]$ satisfies
	\[
\Phi(u) = \begin{cases} &\exp\left(- \ga^\al |u|^{\al} \left[ 1 + i b \tan\left(\frac{\pi \al}{2} \right) \text{sign}(u)  \left(|\ga u|^{1-\al} - 1 \right) \right]  + i u \de \right), ~ \al \neq 1,\\ &\exp\left(- \ga |u| \left[ 1 + i b \frac 2 \pi \text{sign}(u) \log(\ga |u|) \right] + iu \de \right), ~ \al =1, \end{cases}
	\]
	where $0 < \al \le 2$ is the stability index, $-1 \le b \le 1$ is a skewness parameter, $\ga > 0$ is the scaling parameter and $\de \in \R$ is used for location.
	\item[(ii)] A L\'evy process $X=(X_t)_{t \ge 0}$ is called \emph{stable} if the distribution of $X_1$ is stable. 
	\item[(iii)] A L\'evy process $X=(X_t)_{t \ge 0}$ is called a \emph{subordinator} if its sample paths are non-decreasing. 
\end{itemize}
\end{definition}

Below we will work with a pure jump stable subordinator $X=(X_t)_{t \ge 0}$, i.e.\ with a driftless subordinator for which $X_1$ has a one-sided stable distribution. The distribution of such a process is determined by two parameters of $X_1$ only (see e.g.\ Example 4.2 in \cite{contan2004}), and in terms of Definition \ref{defstab} this corresponds to setting $b = 1$ and $\de = 0$, while $\ga > 0$ remains a flexible scaling parameter and $0 < \al < 1$ necessarily because the jumps need to be summable.

\begin{lemma} \label{lemPn}
Let $Z=(Z_t)_{t \ge 0}$ be a  pure jump $\al$-stable subordinator and \(P_n = \sum_{i=1}^n \mathbbm{1}_{\{ |\Delta_i^n Z| > K_n \}}\) as above. Then:
\begin{itemize}
	\item[(i)] $\P(P_n = 0) \longrightarrow 0$ as $n \to \infty$.
	\item[(ii)] For 
	\[
	N_n = C_{\al, \ga} n^{\alpha\varpi}, ~ C_{\al, \ga} =  \ga^\al \frac 2{\pi}\sin\left(\frac{\pi \al}{2} \right) \Ga(\al),
	\]
	we have, as $n \to \infty$,
  \[
    \frac{P_n}{N_n} \pn 1.
  \]
\end{itemize}
\end{lemma}

\begin{theorem} \label{thm1}
Let $Z=(Z_t)_{t \ge 0}$ be a pure jump $\al$-stable subordinator with $0 < \al < 1$, and assume that $\sup_{0 \le s \le t \le 1} \E\left[\left|\si_t - \si_s\right|^p\right|] \le C_p (t-s)^{\ga p}$ for some fixed $\ga > 0$, any $p > 0$ and a constant $C_p > 0$. 
\begin{itemize}
	\item[(i)] Let $\varpi < \frac{1}{2 - \alpha}.$ Then 
	\[
	\hat \be_n - \be = O_\P\left(\max\left(n^{\frac{(2-\al)\varpi - 1}2}, n^{-\frac{\al \varpi}2}\right)\right)
	\]
	and hence $\hat \be_n - \be \pn 0$. 
	\item[(ii)] Suppose that $\frac{1}{2} < \varpi < \frac{1}{2 - \alpha}.$ Then we have the stable convergence 
	\[
	n^{\frac{1-(2-\al)\varpi}2} \left(\hat \be_n - \be \right) \tols \sqrt{\frac{\al}{(\alpha+2)C_{\ga, \al}}\int_{0}^{1} \si_s^2 ds} \cdot Y
	\]
	where $Y$ is a standard normal random variable defined on an appropriate extension of the probability space and independent of $\mathscr F$. 
\end{itemize}
\end{theorem}

\begin{remark} \label{rem1}
\begin{itemize}
	\item[(i)] In its current form Theorem \ref{thm1} is infeasible as it contains several unknown quantities, but estimation of all of these is possible in principle.  Regarding $\ga$ and $\al$ one can rely on the estimation techniques discussed in \cite{nolan2020}, e.g., working with maximum-likelihood methods or using estimation via characteristic functions. Ultimately, the $\De_i^n Z$ are i.i.d.\ stable variables and as such allow for a consistent estimation of the underlying parameters. Also, estimation of the integrated volatility $\int_0^1 \si_s^2 ds$ can be based on observations of $\De_i^n Y$, using the thresholding approach from \cite{mancini2009}. 
	\item[(ii)] The choice of $\varpi$ in Theorem \ref{thm1}(ii) is atypical, as threshold models usually rely on $\varpi < \frac 12$ in order for jump parts to dominate over increments of a Brownian martingale. Such a choice was in principle possible as well, and in this case the leading term in a decomposition of $\hat \be_n - \be$ would be purely governed by the jumps in $Y$ and $Z$. If we compare with (\ref{decomps}), a key step in the derivation of a central limit theorem then was to derive a limit of 
	\[
	\frac{\si_\eta^2}{N_n^2} \sum_{i=1}^n \E \left[ \frac{\sum_{(i-1)\De_n < s \le i\De_n} |\Delta Z_s|^2}{|\Delta_i^n Z|^2}  \mathbbm{1}_{\{ |\Delta_i^n Z| > K_n \}}  \right].
	\]
To the best of our knowledge, the asymptotic behaviour of this term is unknown, and it is even unclear how it would behave in general. We therefore choose $Z$ to be a subordinator, allowing for a simple upper bound for this term, and then choose $\varpi > \frac 12$ in order for this term to be of an asymptotically smaller order. 
	\item[(iii)] The assumption of a strictly stable L\'evy process $Z$ is restrictive, too, but is at least to a certain extent reasonable from a financial perspective; see e.g.\ the models discussed in the works of \cite{schoutens2003} and \cite{contan2004}. Methodologically, we need it for the asymptotic behaviour of 
\[
\E \left[\left| \Delta_i^n Z\right|^{-p} \mathbbm{1}_{\{ |\Delta_i^n Z| > K_n \}}  \right],
\]
see (\ref{lct}), which is again difficult to derive without additional assumptions. 
\end{itemize}

\end{remark}

\subsection{The infinite time horizon}
Ultimately, by working with an infinite time horizon one can be more flexible regarding the nature of the regressor $Z$, as jumps can be observed over the entire positive half line. We hence allow for a setting which mirrors the model from \cite{lietal2017a} and is closer to standard real world applications. Of course, the statistical model from Section \ref{secfin} remains possible as well. 

Hence, let $(Y, Z)$ be as in (\ref{defYZ}), more specifically
\begin{align*}
 &Y_t=Y_0+ \int_0^t a_s ds+\int_0^t \sigma_s dW_s+ \sum_{s \le t} (\beta+\eta_s) \Delta Z_s, \nonumber \\ 
  &Z_t = Z_0+ \int_0^t a'_s ds+\int_0^t \sigma'_s dW'_s+J_t,\quad  t \ge 0,
\end{align*}
where $a, a', \si, \si'$ are all assumed to be  bounded and c\`adl\`ag, and we assume \(\eta\) still to be i.i.d.\ bounded with zero mean and variance \(\sigma^2_\eta\), independent of all other processes. 

What is structurally different, though, is the assumptions on the jump process $J_t$. Here we allow it to be compound Poisson with a time-varying intensity, i.e.\ we have 
\[
J_t = \sum_{i=1}^{R_t}\xi_i
\]
with the following conditions in place. 

\begin{condition} \label{condFA}
\begin{itemize}
	\item[(i)] The process $R=(R_t)_{t \ge 0}$ is an inhomogeneous Poisson process with a non-negative deterministic intensity process $\la=(\la_t)_{t \ge 0}$ which satisfies 
\begin{equation*}
  \limsup_{n\to\infty} \sup_{1 \le i \le n} \frac{1}{\Delta_n} \int_{(i-1)\De_n}^{i \De_n} \lambda_s ds \le \lambda^* 
\end{equation*}
for some $\la^* > 0$ as well as for which there exists some $\overline{\lambda} > 0$ such that
\begin{equation*}
  \frac{1}{x} \int_0^{x} \lambda_s ds \longrightarrow \overline{\lambda} \quad \text{as } x \to \infty.
\end{equation*}
	\item[(ii)] The \(\xi_i\) are i.i.d., independent of all other processes involved, with a three times differentiable, even density \(f\) that satisfies \(f(0)=0\). We further assume that \(\xi\) has a finite $(2+\de)$th moment, $\de > 0$, and that the $f^{(j)}$ are bounded and integrable for all $0 \le j \le 3$.
\end{itemize}
\end{condition}

Regarding the sampling scheme we assume to observe $n$ high-frequency increments with an infinite time horizon, i.e.\ we have \(\Delta_n \longrightarrow 0\) and \(T_n = n\Delta_n \longrightarrow \infty\), with the additional 
  \(n\Delta_n^2 = T_n \Delta_n \longrightarrow 0\)
in place. We further set \(K_n\),  \(P_n\) and $\hat{\beta}_n$ as above. 

\begin{lemma} \label{lemRn}
Let \(K_n =\Delta_n^\varpi \) for some $\frac 13<\varpi < \frac 12$.
\begin{itemize}
	\item[(i)] There exists some $C > 0$ such that 
	\[
	\sup_{1 \le i \le n} \E \left[\left| \mathbbm{1}_{\{ |\Delta_i^n Z| > K_n \}} - \mathbbm{1}_{\{ \Delta_i^n R = 1\}}\right|  \right] \le C \De_n^2
	\]
	for all $n$ large enough. 
	\item[(ii)] For $N_n = n \De_n \overline \la$ we have, as $n \to \infty$,
  \[
    \frac{P_n}{N_n} \pn 1.
  \]
\end{itemize}
\end{lemma}

\begin{theorem}
\label{thm2}
Under the assumptions from Condition \ref{condFA} and with $\frac 13<\varpi < \frac 12$ we have
\begin{equation*}
  \sqrt{N_n} (\hat{\beta}_n -\beta) \tol  \mathcal{N}\left(0, \sigma_\eta^2\right).
\end{equation*}
\end{theorem}

\begin{remark}
\begin{itemize}
	\item[(i)] In the infinite time horizon setup both the central limit theorem itself and the condition of $\varpi$ look much more natural and standard compared with Theorem \ref{thm1} and Remark \ref{rem1}. In fact, we can choose the standard $\varpi < \frac 12$ as we now want to only have those increments involved which are essentially unaffected by the continuous martingale part. As discussed e.g.\ in \cite{figman2019}, such a choice of the threshold allows for an almost surely correct identification of the jump times, and the corresponding increments converge to the correct jump sizes as well. 
	\item[(ii)] Following the above discussion it is very plausible that $\hat \be_n$ indeed reaches the optimal variance in this model. Suppose that the jump sizes and jump times of $(Y,Z)$ could be observed exactly, i.e.\ we have observations $(\De Y_{\tau_k}, \De Z_{\tau_k})$ with
	\[
	\De Y_{\tau_k} = \be \De Z_{\tau_k} + \eta_{\tau_k} \De Z_{\tau_k},
	\]
where the $\tau_1, \ldots, \tau_{R_{T_n}}$ denote the true jump times until $T_n.$ In this case we are dealing with a linear regression involving heteroscedastic errors, and specifically $\Var\left(\eta_{\tau_k} \De Z_{\tau_k}\middle \vert \De Z_{\tau_k} \right) = \si^2_\eta \left(  \De Z_{\tau_k}\right)^2.$ Following e.g.\ Appendix C.11 in \cite{montetal2001} it is clear that the weighted least squares estimator, with weights proportional to the inverse standard deviation, is the best linear unbiased estimator. In this case, as the $\De Z_{\tau_k}$ are known, it becomes 
\[
\widetilde \be_n = \frac 1{R_{T_n}} \sum_{k=1}^{R_{T_n}} \frac{\De Y_{\tau_k}}{\De Z_{\tau_k}} = \be + \frac 1{R_{T_n}} \sum_{k=1}^{R_{T_n}} \eta_{\tau_k}.
\]
Clearly, 
\[
\sqrt{R_{T_n}} \left(\widetilde \be_n - \be \right) \tol \nN\left(0, \si_\eta^2 \right).
\]
With the aid of Lemma \ref{lemRn} it is easy to derive $\frac{R_{T_n}}{N_n} \pn 1$, and hence Theorem \ref{thm1} shows that $\hat \be_n$ shares the asymptotic behaviour with the ``oracle'' estimator $\widetilde \be_n$. 
\item[(iii)] This statistical behaviour also explains why we work with a simpler model than (\ref{defju}) in the sense that we do not consider additional idiosyncratic jumps $\De X_t$. For a large $n$, if we only restrict to those intervals with $\De_i^n Z > K_n$, we identify the jump points of $Z$ correctly, and as $X$ jumps independently of $Z$, these intervals asymptotically are not affected by jumps of $X$. A similar feature holds in \cite{lietal2017a} as well. 
\end{itemize}
\end{remark}

\section{Simulation study} \label{sec:sim}
\def\theequation{4.\arabic{equation}}
\setcounter{equation}{0}

Throughout the simulation study we work with a dependent process of the form
\begin{align*}
  Y_t = Y_0 + a t + \sigma W_t + \sum_{0 < s \le t} (\beta + \eta_s) \De Z_s,
\end{align*}
where the noise variables \(\eta_t\) are drawn from  \(\mathcal{U}(-0.5, 0.5)\). We also set $b = 1.5$, $\si = 10$, $\be = 2$, and we conduct 5000 Monte Carlo replications and choose $n=4000$ in all cases. 

Let us start with a discussion of the behaviour of $\hat \be_n$ in the situation of Theorem \ref{thm1}. For this model we choose $n^{-1} = \De_n$, and the process $Z$ is assumed to be an $\al$-stable subordinator with $\ga = 1$. For its implementation we follow \cite{chaetal1976} and simulate i.i.d.\ stable variables from which we obtain increments $\De_i^n Z$ via the self-similarity property. 

As discussed in Remark \ref{rem1}(ii) the central limit theorem is obtained by choosing $\varpi > \frac 12$ as in this case an interplay of the Brownian increment and the stable process gives the leading term asymptotically. On the other hand, the condition $\varpi < \frac{1}{2-\al}$ forces the parameter to be chosen only slightly larger than $\frac 12$, in particular for a small $\al$, and in this case the contribution of the smaller order terms might not be negligible in practice.

Below we show the results in the two cases $\al=0.3$ and $\al=0.7$, for various values of $\frac 12 < \varpi < \frac 1{2-\al}$. Specifically, we choose three values of $\varpi$: Two of them very close to the above boundaries and one exactly in the middle of the two. 

\begin{table}[h] 
\centering
  \begin{tabular}{|c |c |c|}
    \hline
    \multicolumn{3}{|c|}{$\beta=2, \alpha=0.3 $} \\
    \hline
    \(\varpi\)& \(\hat{\beta}\) & MSE \\
    \hline
    0.501 & 1.9471 & 5.5177 \\
    0.544 & 2.0885  & 10.0971 \\
    0.587 & 2.0147  & 16.8985 \\
    \hline
  \end{tabular}
	\qquad 
	 \begin{tabular}{|c |c |c|}
    \hline
    \multicolumn{3}{|c|}{$\beta=2, \alpha=0.7$} \\
    \hline
    \(\varpi\)& \(\hat{\beta}\) & MSE \\
    \hline
    0.501 & 1.9563 & 2.1436 \\
    0.634 & 2.1535 & 8.5536 \\
    0.768 & 1.8571 & 35.1001 \\
    \hline
  \end{tabular}
	\caption{Empirical mean of $\hat \be_n$ and mean squared error of $\hat \be_n - \be$ in the framework of Theorem \ref{thm1} for $\al = 0.3$ (left) and $\al = 0.7$ (right), $n=4000$ and various values of $\varpi$}
	 \label{t1}
\end{table}

\begin{figure}[!t] 
  \centering
  \includegraphics[width=0.9\textwidth]{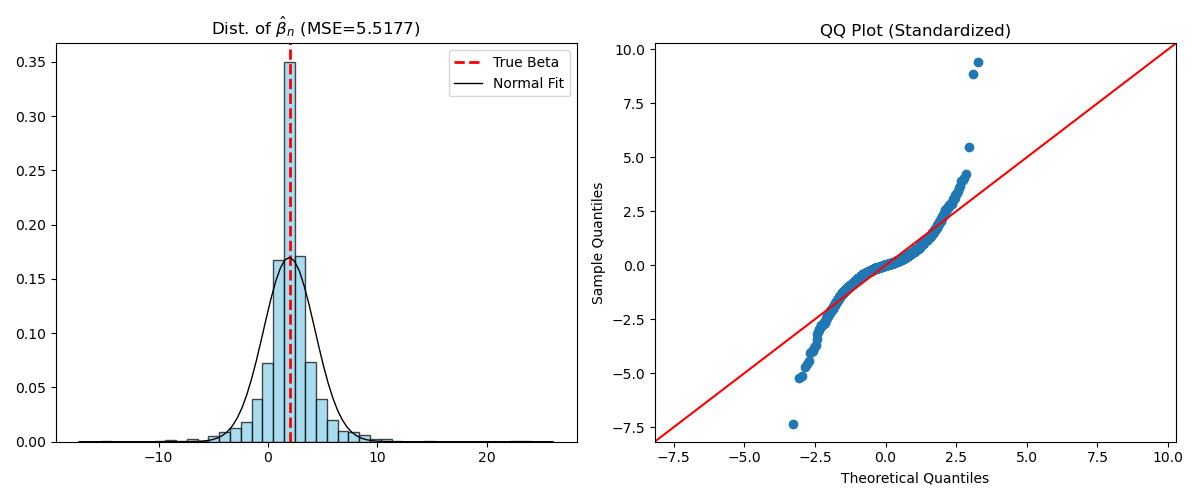}
  \includegraphics[width=0.9\textwidth]{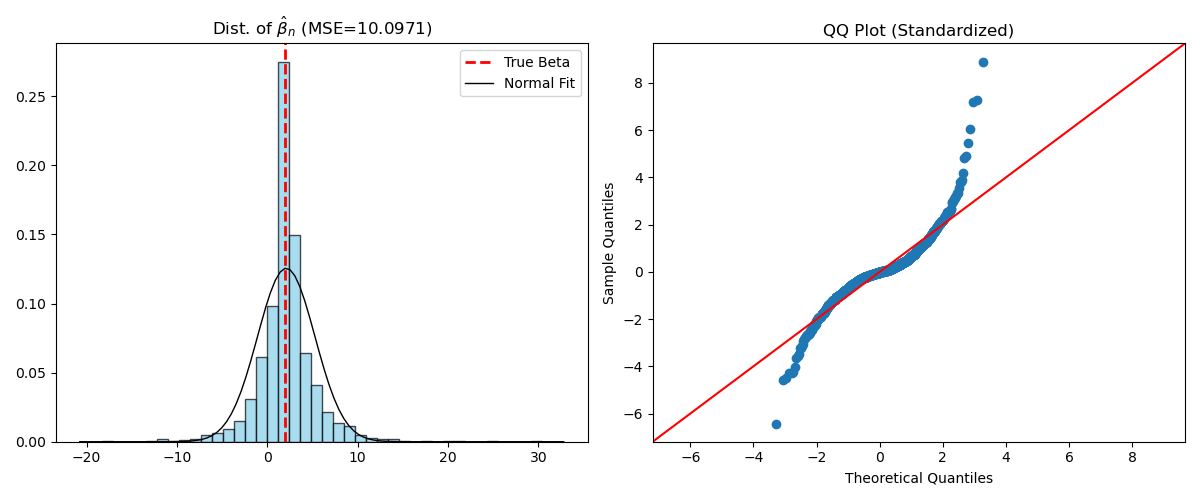}
  \includegraphics[width=0.9\textwidth]{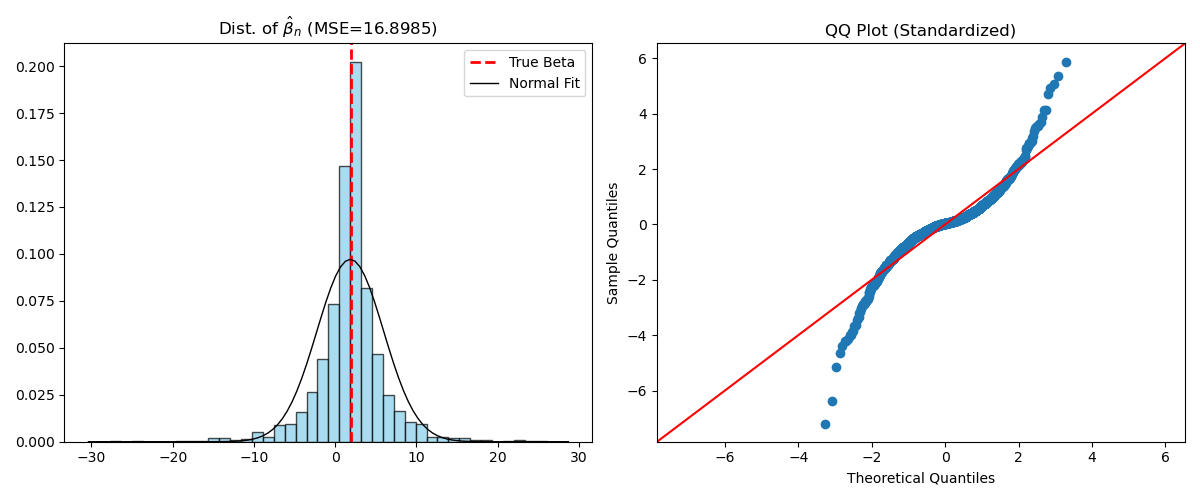}
  \caption{Normal approximation and qq plot for \(\alpha=0.3\). From top to bottom: \(\varpi = 0.501\), \(\varpi = 0.544\), \( \varpi = 0.587\).}
	\label{f1}
\end{figure}

It is visible from Table \ref{t1} that $\hat \be_n$ is in principle a reasonable estimator for $\be$, and in all cases. This is particularly relevant as $\varpi$ is in general unknown and hence a choice of $\varpi$ close to the lower boundary is of importance. The rather large mean squared error in both tables is due to the fact that the asymptotic variance is proportional to $\si^2$ which is of a relevant size due to our choice of $\si=10$.  Note also that the increase in the mean squared error from $\varpi_1$ to the next larger $\varpi_2$ is to be expected as $\varpi$ enters the convergence rate in Theorem \ref{thm1} explicitly. In particular, the increase in the variance should roughly correspond to a factor $n^{(2-\al)(\varpi_2-\varpi_1)}$ which for our choice of $n = 4000$ is 1.834 (for $\al = 0.3$) and 4.196 (for $\al =0.7$). This phenomenon is well reflected in Table \ref{t1}.

\begin{figure}[!t]
  \centering
  \includegraphics[width=0.9\textwidth]{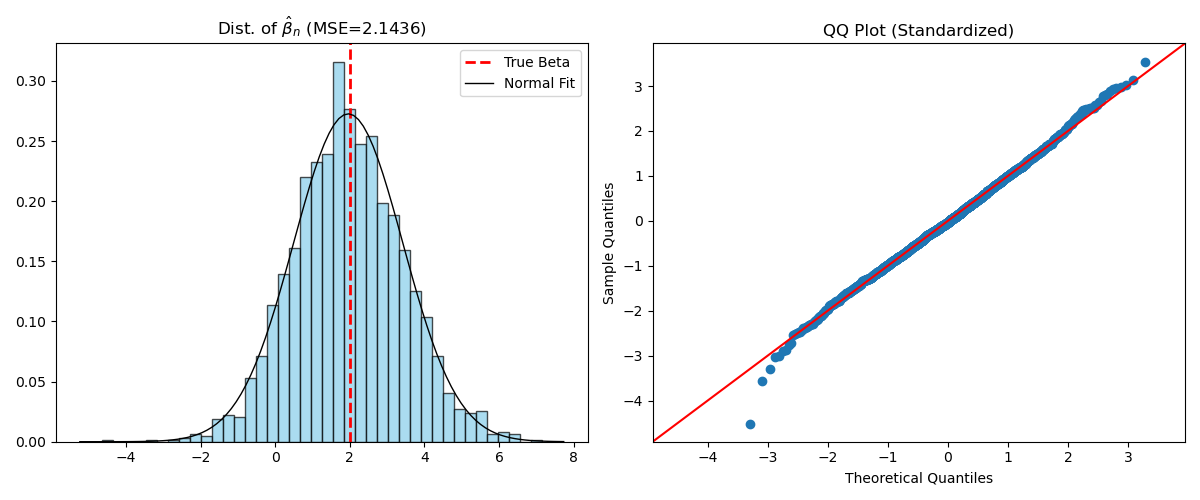}
  \includegraphics[width=0.9\textwidth]{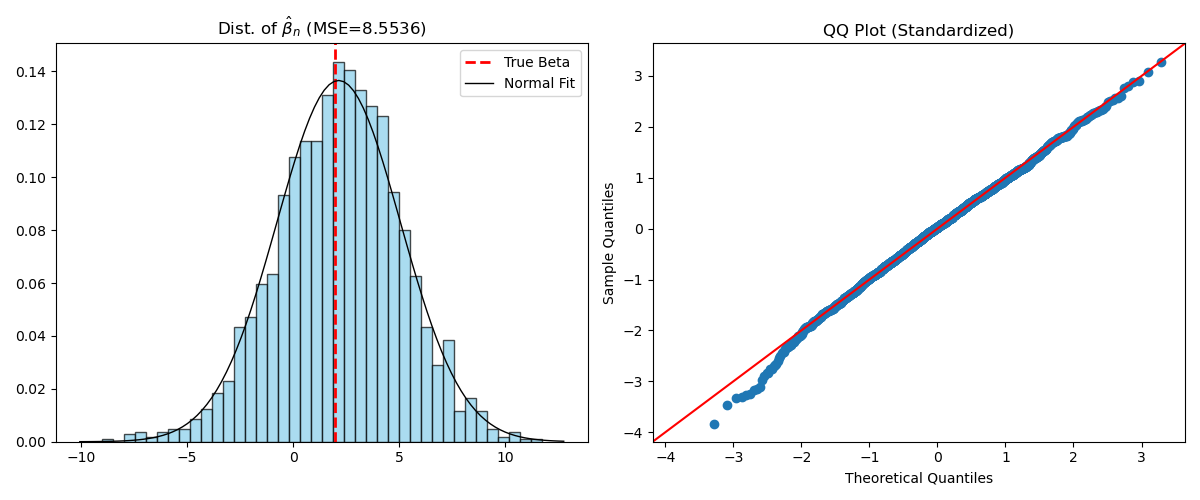}
  \includegraphics[width=0.9\textwidth]{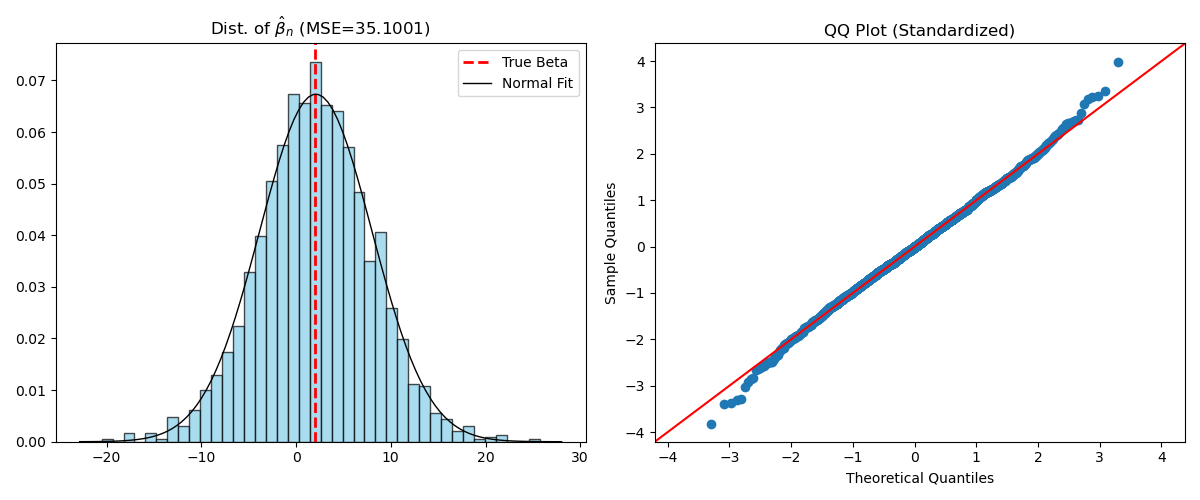}
  \caption{Normal approximation and qq plot for \(\alpha=0.7\). From top to bottom: \(\varpi = 0.501\), \(\varpi = 0.634\), \( \varpi = 0.768\).}
	 \label{f2}
\end{figure}

The fit of the normal distribution is assessed in Figure \ref{f1} for $\al = 0.3$ and in Figure \ref{f2} for $\al = 0.7$, and we see a markedly different behaviour in the two cases. For $\al=0.3$ note that only very few increments of $Z$ actually exceed the threshold $K_n$. Comparing with Lemma \ref{lemPn} the expected number $N_n$ grows as $n^{\al \varpi}$ and takes values in the single digits for our choices of $\al, \ga$ and $n$. Once $\al$ becomes sizeably larger and enough increments contribute to the average, the normal approximation begins to work fairly well. 

In the infinite time horizon case we work with a finite activity jump process $Z$, and we set
\[
  Z_t = Z_0 + a' t + \sigma' W'_t +  J_t
\]
with $a' = 1$, $\si' = 0.5$ and where the Brownian motion $W'$ is independent of $W$. Here, $J_t$ is a compound Poisson process with varying $\la$ and the jump size density is the one of a symmetric beta distribution with parameter 4, and we choose $\De_n$ such that $n \De_n^{\frac 32} = 1$ in order to ensure $n \De_n^2 \longrightarrow 0$ and $T_n = n\De_n \longrightarrow \infty$. 

\begin{table}[h] 
  \centering
  \begin{tabular}{| c| c | c  |}
    \hline
    \multicolumn{3}{|c|}{$\beta=2, \lambda=5$} \\
    \hline
    \(\varpi\)& \(\hat{\beta}\) & MSE \\
    \hline
    0.334 & 1.9981 & 0.0285 \\
    0.416 & 2.0020 & 0.0333 \\
    0.499 & 1.9917 & 0.0387 \\
    \hline
  \end{tabular}
	\qquad
	\begin{tabular}{| c| c | c  |}
    \hline
    \multicolumn{3}{|c|}{$\beta=2, \lambda=100$} \\
    \hline
    \(\varpi\)& \(\hat{\beta}\) & MSE \\
    \hline
    0.334 & 1.9970 & 0.0015 \\
    0.416 & 2.0016 & 0.0018 \\
    0.499 & 2.0047 & 0.0021 \\
    \hline
  \end{tabular}
				\caption{Empirical mean of $\hat \be_n$ and mean squared error of $\hat \be_n - \be$ in the framework of Theorem \ref{thm2} for $\la = 5$ (left), $\la = 100$ (right), $n=4000$ and various values of $\varpi$}
	\label{t3}
\end{table}

It is visible from Table \ref{t3} that the estimation of $\be$ works very well in this situation, almost regardless of the choice of $\varpi$. Even for a moderate choice of $\la = 5$ we expect about 80 intervals to contain jumps, which is large enough for consistency to kick in, and a larger choice of $\la$ then leads to a decrease in the mean squared error as the asymptotic variance is proportional to the inverse of $N_n$.

\begin{figure}[!t]
  \centering
	\includegraphics[width=0.9\textwidth]{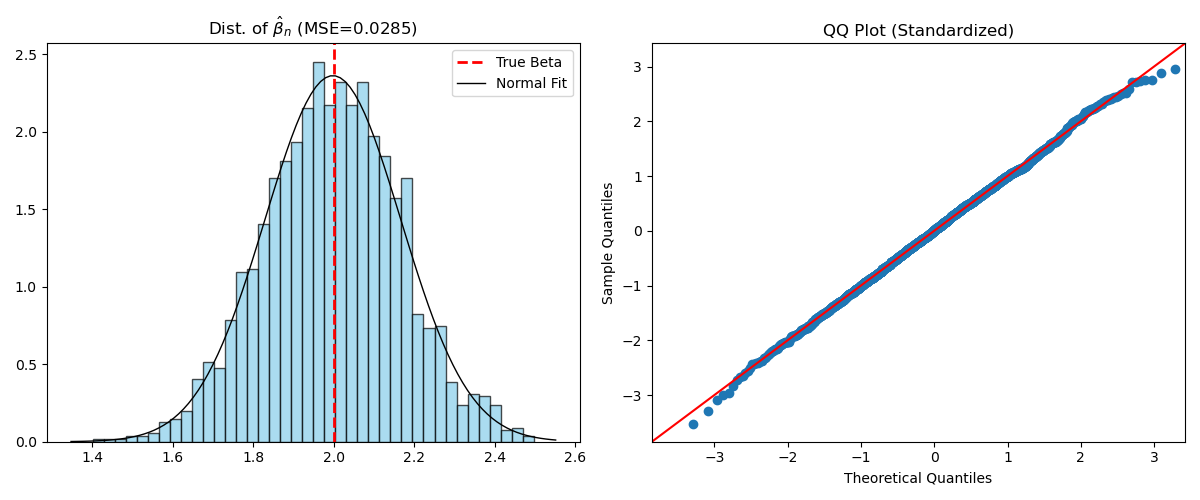}
	\includegraphics[width=0.9\textwidth]{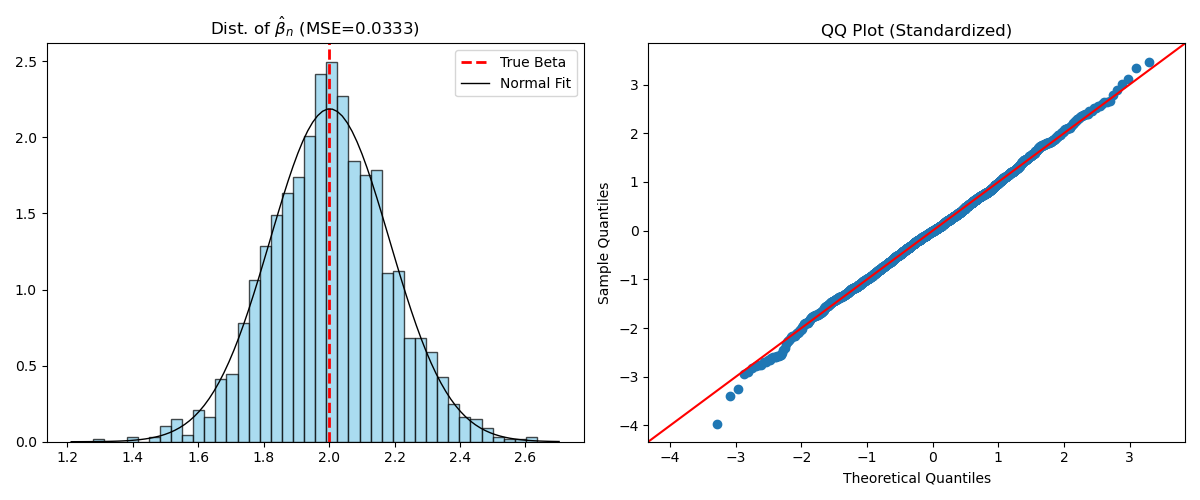}
  \includegraphics[width=0.9\textwidth]{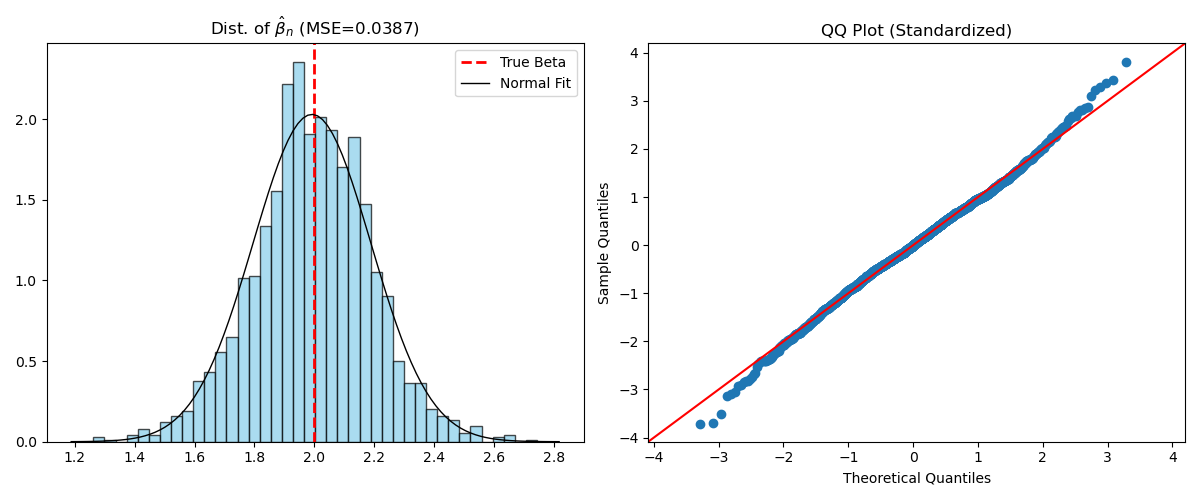}
  \caption{Normal approximation and qq plot for \(\la=5\). From top to bottom: \(\varpi = 0.334\), \(\varpi = 0.416\), \( \varpi = 0.499\).}
	 \label{f3}
\end{figure}

\begin{figure}[!t]
  \centering
  \includegraphics[width=0.9\textwidth]{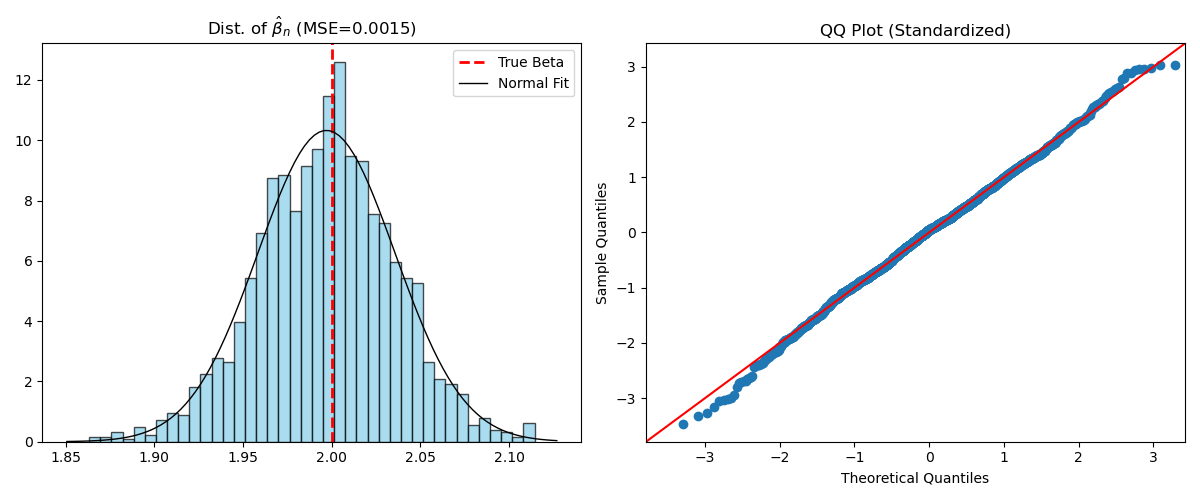}
  \includegraphics[width=0.9\textwidth]{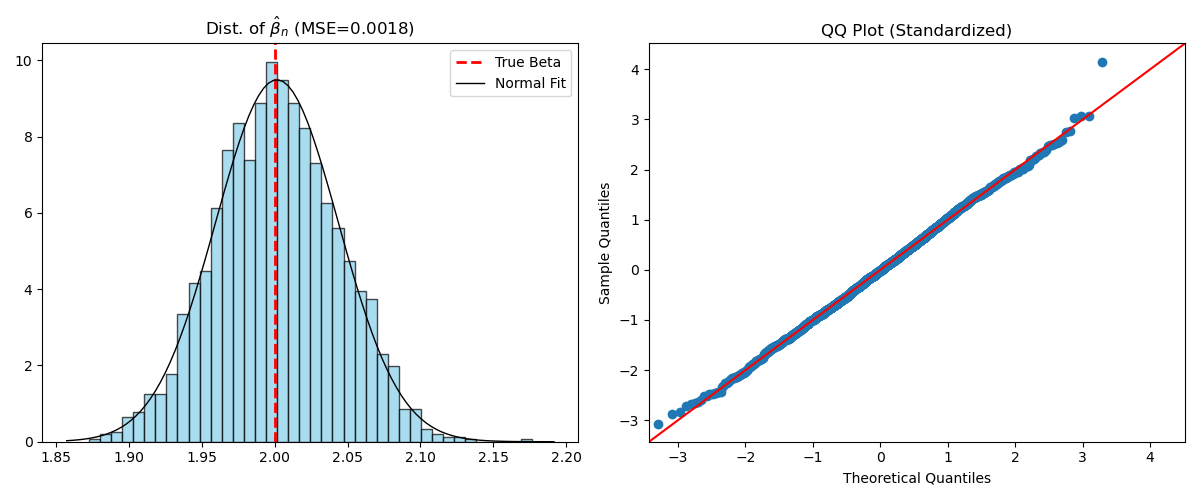}
  \includegraphics[width=0.9\textwidth]{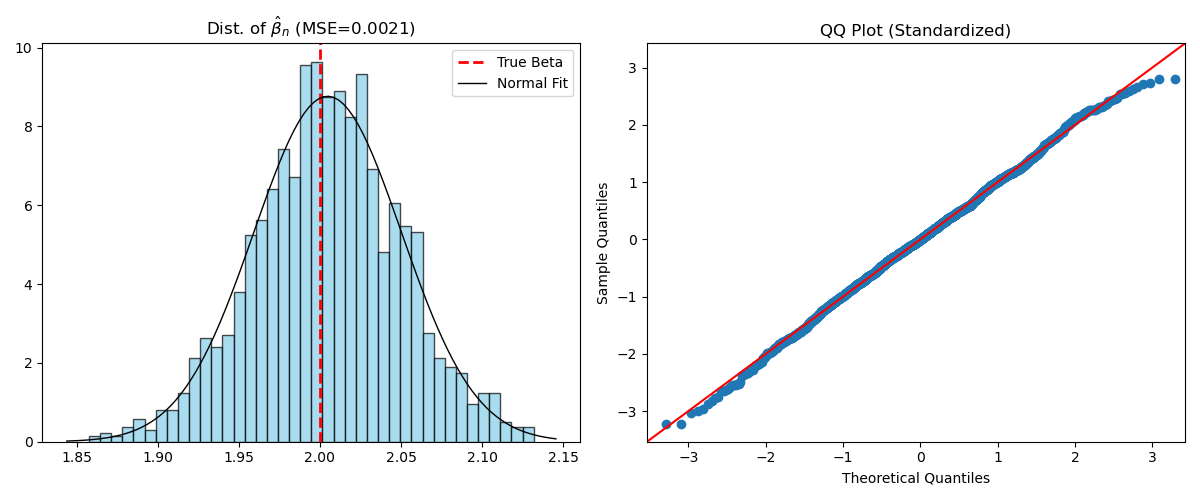}
  \caption{Normal approximation and qq plot for \(\la=100\). From top to bottom: \(\varpi = 0.334\), \(\varpi = 0.416\), \( \varpi = 0.499\).}
	 \label{f4}
\end{figure}

These findings are supported by the corresponding plots on the accuracy of the normal approximation which are given in Figure \ref{f3} and Figure \ref{f4}. It is visible that this approximation works fairly well in both cases, and additional simulation results not reported here show that the statistical properties are not very sensitive to the choice of $\si'$, suggesting that the threshold mechanisms works fine even for larger values of the underlying volatility.

\section{Proofs} \label{sec:proof}
\def\theequation{5.\arabic{equation}}
\setcounter{equation}{0}

\subsection{Proof of Lemma \ref{lemPn}}

The proof heavily relies on the self-similarity property of a stable process, in connection with a tail approximation for stable random variables with $\al < 2$ and $b > -1$. Precisely, Theorem 1.2 in \cite{nolan2020} gives 
\begin{align}
\label{tailapp}
\P(Z_1 > x) \sim \ga^\al c_\al (1+\be) x^{-\alpha}
\end{align}
as \(x \to \infty\), where $c_\al = \frac 1{\pi}\sin\left(\frac{\pi \al}{2} \right) \Ga(\al)$.

Note for a L\'evy process $X=(X_t)_{t \ge 0}$ that its characteristic function is always of the form $\Phi_{X_t}(u) = \exp(t \kappa(u))$. Hence, comparing with Definition \ref{defstab}(i) it is clear that each $Z_t$, $t \ge 0$, is stable as well, with parameters $\al, b, \ga t^{1/\al}, \de t$. In our case of a pure jump $\al$-stable subordinator with $b=1, \de = 0$ only the scaling parameter changes, and hence it is easy to deduce that  
\(Z_{t} \stackrel{\mathcal{L}}{=} t^{1/\alpha} Z_1\).

Now, let \(p_n = \P(\Delta_i^n Z > K_n)\). By the stationarity of increments and from the preparation above we get the self-similarity property \(\De_i^n Z \stackrel{\mathcal{L}}{=} \De_n^{1/\alpha} Z_1 = n^{-1/\al} Z_1\), hence
  \begin{align*}
    p_n &= \P\left( \Delta_1^n Z > K_n \right) \\
    &= \P\left( n^{-1/\alpha} Z_1 > K_n \right) = \P\left( Z_1 > n^{1/\alpha} K_n \right) = \P\left( Z_1 > n^{1/\alpha-\varpi} \right)
  \end{align*}
	by definition of $K_n$. From $\al < 1$ and $\varpi < 1$ it obvious that $n^{1/\alpha-\varpi} \longrightarrow \infty$. (\ref{tailapp}) then gives 
  \begin{align*}
    p_n &\sim  2 \ga^\al c_\al \left( n^{\frac{1}{\alpha} - \varpi} \right)^{-\alpha} = C_{\al, \ga} \left( n^{\frac{1}{\alpha} - \varpi} \right)^{-\alpha} = C_{\al, \ga} n^{\alpha\varpi-1}.
  \end{align*}
  We also have 
  \begin{align*}
    \E[P_n] &= \sum_{i=1}^n p_n = n p_n \sim n  C_{\al, \ga} n^{ \alpha\varpi -1} = C_{\al, \ga} n^{\alpha\vpi} = N_n, 
  \end{align*}
  implying \(N_n \sim \E[P_n] \longrightarrow \infty\) as \(n \longrightarrow \infty\).

 To prove (i) we bound the probability \(\P(P_n = 0)\), using the elementary inequality \(\log(1 - x) \le -x\), which is valid for \(x<1\), via
  \begin{align*}
    \P(P_n = 0) &= (1 - p_n)^n = \exp(n \log(1-p_n)) \le \exp( -n p_n ).
  \end{align*}
  Using \(n p_n \sim C_{\ga, \la} n^{\alpha\vpi} \longrightarrow \infty\) we then get
  \begin{align*}
    \exp( -n p_n ) \longrightarrow 0 \quad \text{as } n \longrightarrow \infty.
  \end{align*}

  To prove (ii) note first that \(P_n \sim \text{Bin}(n, p_n)\). We apply Chebyshev's inequality to the ratio \(\frac{P_n}{\E[P_n]}\). First,
  \begin{align*}
    \frac{\text{Var}(P_n)}{(\E[P_n])^2} &= \frac{n p_n (1 - p_n)}{(n p_n)^2} = \frac{1 - p_n}{n p_n}.
  \end{align*}
  From \(n p_n \longrightarrow \infty\) we get
  \begin{align*}
    \lim_{n \to \infty} \frac{\Var(P_n)}{(\E[P_n])^2} = 0.
  \end{align*}
  By Chebyshev's inequality, for any \(\epsilon > 0\):
  \begin{align*}
    \P\left( \left| \frac{P_n}{\E[P_n]} - 1 \right | > \epsilon \right) \le \frac{1}{\epsilon^2} \frac{\text{Var}(P_n)}{(\E[P_n])^2} \longrightarrow 0.
  \end{align*}
  The claim follows since $N_n \sim \E[P_n]$.

	\subsection{Proof of Theorem \ref{thm1}}
	By a standard localisation procedure as in Section 4.4.1 of \cite{jacpro12} we can assume boundedness of the processes $a$ and $\si$, and we use $C$ throughout the proof to denote a universal constant. 
	
	Now with the notation \(\Omega_n = \{P_n \ge 1\}\) we decompose the error as follows:
  \begin{align*}
    \hat{\beta}_n - \beta &= \frac{1}{P_n} \sum_{i=1}^n \left( \frac{\int_{(i-1)\De_n}^{i \De_n} a_s ds}{\Delta_i^n Z} + \frac{\int_{(i-1)\De_n}^{i \De_n} \si_s dW_s}{\Delta_i^n Z} \right. \\ &~~~~~~~~~~~~~~~~~~~~~~~~~+ \left. \frac{\sum_{(i-1)\De_n < s \le i\De_n} \eta_s \Delta Z_s}{\Delta_i^n Z} \right) \mathbbm{1}_{\{ |\Delta_i^n Z| > K_n \}}  \mathbbm{1}_{ \Omega_n } + (0-\beta)\mathbbm{1}_{ \Omega_n^\mathsf{c} }.
  \end{align*}
	From Lemma \ref{lemPn} we know that $\P(\Om_n^{\mathsf{c}}) \longrightarrow 0$ as well as $\frac{P_n}{N_n} \pn 0$. Hence, from Slutsky's lemma, it is sufficient to discuss the three summands in  
	\[
	 \frac{1}{P_n} \sum_{i=1}^n \left( \frac{\int_{(i-1)\De_n}^{i \De_n} a_s ds}{\Delta_i^n Z} + \frac{\int_{(i-1)\De_n}^{i \De_n} \si_s dW_s}{\Delta_i^n Z} + \frac{\sum_{(i-1)\De_n < s \le i\De_n} \eta_s \Delta Z_s}{\Delta_i^n Z} \right) \mathbbm{1}_{\{ |\Delta_i^n Z| > K_n \}}  \mathbbm{1}_{ \Omega_n }
	\]
	separately. We also may replace $P_n$ by $N_n$ whenever reasonable, and if we do so the indicator over $\Om_n$ may be dropped as each $\{|\Delta_i^n Z| > K_n\}$ is a subset of $\Om_n$.

For first term involving the drift, note that 
\begin{align*}
\frac{1}{P_n} \sum_{i=1}^n \left| \frac{\int_{(i-1)\De_n}^{i \De_n} a_s ds}{\Delta_i^n Z} \right| \mathbbm{1}_{\{ |\Delta_i^n Z| > K_n \}} \mathbbm{1}_{ \Omega_n }  \le \frac{1}{P_n} \sum_{i=1}^n  \frac{C\De_n}{K_n} \mathbbm{1}_{\{ |\Delta_i^n Z| > K_n \}} \mathbbm{1}_{ \Omega_n }  = C n^{\varpi -1} \mathbbm{1}_{ \Omega_n }
\end{align*}
by definition of $P_n$. 

In order to discuss the term involving the continuous martingale part, let us first discuss a general formula for
\[
\E \left[\left| \Delta_i^n Z\right|^{-p} \mathbbm{1}_{\{ |\Delta_i^n Z| > K_n \}}  \right] = n^{p/\al}\E \left[\left| Z_1\right|^{-p} \mathbbm{1}_{\{ |Z_1| > n^{1/\al-\varpi} \}}  \right], \quad p \ge 1,
\] 
where the identity stems from the self-similarity property. We here use Theorem 1.2 in \cite{nolan2020} again. For the tail density of a stable random variable we have 
\begin{align}
\label{tailapp2}
f(x) \sim \al \ga^\al c_\al (1+b) x^{-(\alpha+1)} 
\end{align}
as \(x \longrightarrow \infty\). The limit comparison test then easily gives
  \begin{align} \label{lct}
   \nonumber \E \left[\left| \Delta_i^n Z\right|^{-p} \mathbbm{1}_{\{ |\Delta_i^n Z| > K_n \}}  \right] &= n^{p/\al}\E \left[\left| Z_1\right|^{-p} \mathbbm{1}_{\{ |Z_1| > n^{1/\al-\varpi} \}}  \right] \\ 
	 \nonumber &\sim   n^{p/\al} 	\int_{n^{1/\al-\varpi}}^{\infty} x^{-p} \al C_{\ga, \al} x^{-(\alpha+1)}dx  \\ &= n^{p/\al} \frac{\al C_{\ga, \al}}{\alpha+p} \left[ -x^{-(\alpha+p)} \right]_{n^{1/\al-\varpi}}^{\infty} = \frac{\al C_{\ga, \al}}{\alpha+p} n^{(\al+p)\varpi - 1}
  \end{align}
where we have used $b = 1$ and the definition of $C_{\ga, \al}$. Note that $\frac 1\al > 1 > \varpi$ easily gives $n^{1/\al-\varpi} \longrightarrow \infty$. 

We now follow the discussion in the proof of Theorem 2.3 in \cite{amoetal2024}, in combination with Lemma 3(ii) and Theorem 6 in \cite{cohen2013}: The martingale representation theorem for jump measures proves that 
\[
 \left|\Delta_i^n Z\right|^{-1} \mathbbm{1}_{\{ |\Delta_i^n Z| > K_n \}}  - \E \left[\left| \Delta_i^n Z\right|^{-1} \mathbbm{1}_{\{ |\Delta_i^n Z| > K_n \}}  \right] = \int_{(i-1)\De_n}^{i\De_n} \int_\R \rho_i^n(x,s) (\mu-\nu)(dx, ds)
\]
where  $\nu$ is the compensator of the random measure $\mu$ associated with the pure jump L\'evy process $Z$ and $\rho_i^n(x,s)$ is square integrable with respect to $\nu$. As the covariation between any continuous and any pure jump martingale is zero we easily get 
\begin{align} \label{exp0}
\E \left[\int_{(i-1)\De_n}^{i\De_n} \int_\R \rho_i^n(x,s) (\mu-\nu)(dx, ds) \int_{(i-1)\De_n}^{i \De_n} \sigma_s dW_s \right] = 0
\end{align}
and we can then conclude that 
\[
\frac{1}{N_n} \sum_{i=1}^n  \frac{\int_{(i-1)\De_n}^{i \De_n} \si_s dW_s}{\Delta_i^n Z} \mathbbm{1}_{\{ |\Delta_i^n Z| > K_n \}}
\]
is a sum of martingale differences with expectation zero, and it is so irrespective of the exact form of $\si$. In particular, we also obtain
\begin{align*}
&\E \left[\left| \frac{1}{N_n} \sum_{i=1}^n  \frac{\int_{(i-1)\De_n}^{i \De_n} \left( \si_s - \si_{(i-1)\De_n}\right) dW_s}{\Delta_i^n Z} \mathbbm{1}_{\{ |\Delta_i^n Z| > K_n \}} \right|^2  \right] \\ =& \frac{1}{N_n^2} \sum_{i=1}^n  \E \left[\left| \frac{\int_{(i-1)\De_n}^{i \De_n} \left(\si_s - \si_{(i-1)\De_n}\right) dW_s}{\Delta_i^n Z} \mathbbm{1}_{\{ |\Delta_i^n Z| > K_n \}}   \right|^2  \right] \\ \le& \frac{1}{N_n^2} \sum_{i=1}^n \E \left[\left| \int_{(i-1)\De_n}^{i \De_n} \left( \si_s - \si_{(i-1)\De_n}\right) dW_s \right|^{2r}  \right]^{\frac 1r} \E \left[\left| \Delta_i^n Z\right|^{-\frac{2r}{r-1}} \mathbbm{1}_{\{ |\Delta_i^n Z| > K_n \}}  \right]^{\frac {r-1}r} \\ \le& C_r \frac{n^{\frac{r-1}r \left( \left(\al + \frac{2r}{r-1}\right)\varpi - 1 \right)}}{N_n^2} \sum_{i=1}^n \E \left[ \left( \int_{(i-1)\De_n}^{i \De_n} \left( \si_s - \si_{(i-1)\De_n}\right)^2 ds \right)^r  \right]^{\frac 1r} \left( 1+ o(1) \right)
\end{align*}
for any integer $r \ge 2$, where we used both the Hölder inequality and the Burkholder-Davis-Gundy inequality for the final two bounds, as well as (\ref{lct}). Several further applications of the Cauchy-Schwarz inequality, together with the uniform bound 
$\E[|\si_s - \si_{(i-1)\De_n}|^p] \le C_p \De_n^{p\ga}$,  prove that the right hand side above is of the order 
\begin{align*}
\frac{n^{\frac{r-1}r \left( \left(\al + \frac{2r}{r-1}\right)\varpi - 1 \right)}}{N_n^2} n \De_n^{1+2 \ga}  
=n^{(2-\al)\varpi -1 + \frac 1r \left(1- \al \varpi \right) - 2 \ga} = o\left(n^{(2-\al)\varpi - 1} \right)
\end{align*}
if we choose $r$ large enough such that $\frac 1r (1-\al \varpi) < 2\ga$.  Hence, 
\begin{align} \label{neg}
\frac{1}{N_n} \sum_{i=1}^n  \frac{\int_{(i-1)\De_n}^{i \De_n} \left( \si_s - \si_{(i-1)\De_n}\right) dW_s}{\Delta_i^n Z} \mathbbm{1}_{\{ |\Delta_i^n Z| > K_n \}}= o_\P\left(n^{\frac{(2-\al)\varpi - 1}2}\right).
\end{align}
%
%
On the other hand, 
\begin{align*}
&\E \left[\left| \frac{1}{N_n} \sum_{i=1}^n  \frac{\si_{(i-1)\De_n} \De_i^n W}{\Delta_i^n Z} \mathbbm{1}_{\{ |\Delta_i^n Z| > K_n \}} \right|^2  \right] = \frac{1}{N_n^2} \sum_{i=1}^n  \E \left[\left| \frac{\si_{(i-1)\De_n} \De_i^n W}{\Delta_i^n Z} \mathbbm{1}_{\{ |\Delta_i^n Z| > K_n \}}   \right|^2  \right],
\end{align*}
and by successive conditioning, and using independence of $\De_i^n W$ and $\De_i^n Z$, 
\begin{align*}
&\E \left[\left| \frac{\si_{(i-1)\De_n} \De_i^n W}{\Delta_i^n Z} \mathbbm{1}_{\{ |\Delta_i^n Z| > K_n \}}   \right|^2  \right]  =  \E \left[ \si_{(i-1)\De_n}^2 \E \left[ \left|\frac{ \De_i^n W}{\Delta_i^n Z} \mathbbm{1}_{\{ |\Delta_i^n Z| > K_n \}}   \right|^2 \middle| \mathscr F_{(i-1)\De_n}  \right] \right]\\ \le & C \E \left[ \left|\De_i^n W\right|^2 \right] \E \left[ \left|\Delta_i^n Z\right|^2 \mathbbm{1}_{\{ |\Delta_i^n Z| > K_n \}} \right]. 
\end{align*}
Together with (\ref{lct}) for $p=2$ this proves 
%
\begin{align*}
&\E \left[\left| \frac{1}{N_n} \sum_{i=1}^n  \frac{\si_{(i-1)\De_n} \De_i^n W}{\Delta_i^n Z} \mathbbm{1}_{\{ |\Delta_i^n Z| > K_n \}} \right|^2  \right] = O\left(\frac{n^{(\al+2)\varpi - 1}}{N_n^2}\right) = O\left(n^{(2-\al)\varpi - 1}\right).
\end{align*}
With a view on (\ref{neg}) it is clear that the latter term is the dominating one in the discussion of the continuous martingale part, and from $n^{2\varpi-2} \le n^{(2-\al)\varpi - 1}$ because of $\al \varpi < 1$ we also see that the latter term above is always of a larger order than the one involving the drift. 

Finally, we have 
\[
\frac{1}{N_n} \sum_{i=1}^n \E \left[\frac{\sum_{(i-1)\De_n < s \le i\De_n} \eta_s \Delta Z_s}{\Delta_i^n Z} \mathbbm{1}_{\{ |\Delta_i^n Z| > K_n \}}  \right]  = 0.
\]
Here we first condition on $ \mathscr{G}_{1}$, the $\si$-field generated by $Z$ until time $t=1$, and then use $\E[\eta_s| \mathscr{G}_{1}] =0$ as $\eta_s$ has mean zero and is independent of $Z$. The same reasoning then gives
\begin{align} \label{decomps} \nonumber
 &\E \left[\left| \frac{1}{N_n} \sum_{i=1}^n \frac{\sum_{(i-1)\De_n < s \le i\De_n} \eta_s \Delta Z_s}{\Delta_i^n Z}  \mathbbm{1}_{\{ |\Delta_i^n Z| > K_n \}} \right|^2  \right] \\ 	&= \frac{\si_\eta^2}{N_n^2} \sum_{i=1}^n \E \left[ \frac{\sum_{(i-1)\De_n < s \le i\De_n} |\Delta Z_s|^2}{|\Delta_i^n Z|^2}  \mathbbm{1}_{\{ |\Delta_i^n Z| > K_n \}}  \right]
\end{align}
Since we assume $Z$ to be a subordinator, however, it is clear that it only has positive jumps. Hence, $\sum_{(i-1)\De_n < s \le i\De_n} |\Delta Z_s|^2 \le |\Delta_i^n Z|^2$ and
\begin{align*}
 &\E \left[\left| \frac{1}{N_n} \sum_{i=1}^n \frac{\sum_{(i-1)\De_n < s \le i\De_n} \eta_s \Delta Z_s}{\Delta_i^n Z}  \mathbbm{1}_{\{ |\Delta_i^n Z| > K_n \}} \right|^2  \right] \\ 	&\le \frac{\si_\eta^2}{N_n^2} \sum_{i=1}^n  \P( |\Delta_i^n Z| > K_n ) \sim \frac{\si_\eta^2}{N_n} = O(n^{-\al \varpi})
\end{align*}
which proves (i).

From the previous calculations and the condition $\frac{1}{2} < \varpi < \frac{1}{2 - \alpha}$ it is clear that 
	\[
	n^{\frac{1-(2-\al)\varpi}2} \left(\hat \be_n - \be \right) = \frac{n^{\frac{1-(2-\al)\varpi}2}}{N_n} \sum_{i=1}^n  \frac{\si_{(i-1)\De_n} \De_i^n W}{\Delta_i^n Z} \mathbbm{1}_{\{ |\Delta_i^n Z| > K_n \}} \left( 1 + o_\P(1)
	\right).
	\]
In order to prove the stable convergence we use Theorem 2.2.15 in \cite{jacpro12} and their notation, and for simplicity we prove (2.2.34), (2.2.36), (2.2.37) and (2.2.39) therein only with $t=1$ as the proof for a general $t \in [0,1]$ works similarly in all cases. We set 
\[
\ze_i^n = \frac{n^{\frac{1-(2-\al)\varpi}2}}{N_n}  \frac{\si_{(i-1)\De_n} \De_i^n W}{\Delta_i^n Z} \mathbbm{1}_{\{ |\Delta_i^n Z| > K_n \}}.
\]
Using independence of $W$ and $Z$ the term above has a vanishing $\mathscr F_{(i-1)\De_n}$-conditional expectation, proving (2.2.34) with $A = 0$. We also have  
\begin{align*}
\sum_{i=1}^n \E\left[ |\ze_i^n|^2 \middle| \mathscr F_{(i-1)\De_n} \right] &=  \frac{n^{{1-(2-\al)\varpi}}}{N_n^2} \sum_{i=1}^n \E \left[\frac{\left|\si_{(i-1)\De_n} \De_i^n W\right|^2}{\left|\Delta_i^n Z\right|^2} \mathbbm{1}_{\{ |\Delta_i^n Z| > K_n \}} \middle| \mathscr F_{(i-1)\De_n}\right] \\ &\sim \frac{n^{{1-(2-\al)\varpi}}}{N_n^2} \sum_{i=1}^n \si_{(i-1)\De_n}^2 \De_n \frac{\al C_{\ga, \al}}{\alpha+2} n^{(\al+2)\varpi - 1}
\end{align*}
by the It\^o isometry and (\ref{lct}). (2.2.36) then follows as the latter quantity converges in probability to $\frac{\al}{(\alpha+2)C_{\ga, \al}}\int_{0}^{1} \si_s^2 ds$, using the assumption on the process $\si$. To prove condition (2.2.37) with $p=3$ we use (\ref{lct}) to obtain 
  \begin{align*}
   \nonumber &\E \left[\left| \Delta_i^n Z\right|^{-3} \mathbbm{1}_{\{ |\Delta_i^n Z| > K_n \}}  \right]\sim  \frac{\al C_{\ga, \al}}{\alpha+3} n^{(\al+3)\varpi - 1}
  \end{align*}
and, hence, 
\begin{align*}
\sum_{i=1}^n \E\left[ |\ze_i^n|^3 \middle| \mathscr F_{(i-1)\De_n} \right] &=  \frac{n^{\frac 32 (1-(2-\al)\varpi)}}{N_n^3} \sum_{i=1}^n \E \left[\frac{\left|\si_{(i-1)\De_n} \De_i^n W\right|^3}{\left|\Delta_i^n Z\right|^3} \mathbbm{1}_{\{ |\Delta_i^n Z| > K_n \}} \middle| \mathscr F_{(i-1)\De_n}\right] \\ &\le C \frac{n^{\frac 32 (1-(2-\al)\varpi)}}{N_n^3} n \De_n^{\frac 32} n^{(\al+3)\varpi - 1} = C n^{-\frac 12 \al \varpi} \longrightarrow 0. 
\end{align*}

We finally prove 
\[
\sum_{i=1}^n \E\left[\ze_i^n \De_i^n M \middle| \mathscr F_{(i-1)\De_n} \right] \pn 0
\]
where $M$ is either $W$ itself or a continuous square integrable martingale orthogonal to $W$. This is sufficient to deduce the remaining condition for stable convergence (2.2.39), e.g.\ according to \cite{jactod2017} and \cite{amoetal2024}. Now, we first have 
\begin{align*}
\left| \E \left[ \frac{\si_{(i-1)\De_n} \left|\De_i^n W\right|^2}{\Delta_i^n Z} \mathbbm{1}_{\{ |\Delta_i^n Z| > K_n \}} \middle| \mathscr F_{(i-1)\De_n} \right] \right| \le \frac{C}{K_n} \De_n \mathbbm{1}_{\{ |\Delta_i^n Z| > K_n \}},
\end{align*}
and then 
\begin{align*}
\left | \sum_{i=1}^n \E\left[\ze_i^n \De_i^n W \middle| \mathscr F_{(i-1)\De_n} \right] \right| \le C \frac{n^{\frac{1-(2-\al)\varpi}2}}{N_n} P_n \frac{\De_n}{K_n} = C\frac{P_n}{N_n} n^{\frac{1-(2-\al)\varpi}2+ \varpi-1} = C\frac{P_n}{N_n} n^{\frac{\al \varpi -1}{2}}.
\end{align*}
The latter converges in probability to zero according to Lemma \ref{lemPn}(ii) and $\al \varpi < 1$. If $M$ is orthogonal to $W$ then $\De_i^n W \De_i^n M$ still defines an increment of a continuous martingale. Hence, using the same reasoning as in (\ref{exp0}) we get 
\begin{align*} 
\E \left[\De_i^n W \De_i^n M \int_{(i-1)\De_n}^{i\De_n} \int_\R \rho_i^n(x,s) (\mu-\nu)(dx, ds) \right] = 0
\end{align*}
It is then simple to deduce 
\[
\E\left[\ze_i^n \De_i^n M \middle| \mathscr F_{(i-1)\De_n} \right] =0. 
\]

\subsection{Proof of Lemma \ref{lemRn}}
To prove (i) let us decompose as follows:
  \begin{align*}
    \left| \mathbbm{1}_{\{ |\Delta_i^n Z| > K_n \}} - \mathbbm{1}_{\{ \Delta_i^n R = 1\}}\right| 
    \le \mathbbm{1}_{\left\{ \Delta_i^n R\ge 2\right\}}
    + \mathbbm{1}_{\left\{ \left|{\Delta_i^n Z}\right| > K_n, \Delta_i^n R = 0 \right\}}
    + \mathbbm{1}_{\left\{ \left|{\Delta_i^n Z}\right| \le K_n , \Delta_i^n R=1\right\}}
  \end{align*}
  and analyse the expectation of each term. 
	
	We use often that the distribution of any $\De_i^n R$ is Poisson with parameter $\De_i^n \La$, setting $\La_t = \int_0^t \la_s ds$. Hence, if we choose $n$ large enough, it is clear by assumption that $\De_i^n \La \le C \De_n$ for some universal $C > 0$, and it is then easy to deduce
	$
	\P(\Delta_i^n R\ge 2) \le C \Delta_n^2$. 
	
	To deal with the second summand above we first condition on $\{ \De_i^n R = 0 \}$ and use the fact that 
	\(\Delta_i^n Z = \int_{(i-1)\De_n}^{i \De_n} a'_s ds + \int_{(i-1)\De_n}^{i \De_n} \sigma'_s dW_s\)
	in this case. Then
      \begin{align*}
        \P(\left|{\Delta_i^n Z}\right| > K_n,  \Delta_i^n R = 0) &= \P(\left| {\Delta_i^n Z} \right| > K_n | \Delta_i^n R = 0)  \P(\Delta_i^n R = 0), 
      \end{align*}
	and the claim then follows easily from
	\begin{align} \label{boundexp}
	 \P(\left| {\Delta_i^n Z} \right| > K_n | \Delta_i^n R = 0) = \P\left(\left|\int_{(i-1)\De_n}^{i \De_n} a'_s ds + \int_{(i-1)\De_n}^{i \De_n} \sigma'_s dW_s \right| > K_n \right) \le C_p \De_n^{p}
	\end{align}
	for any fixed $p > 0$, e.g.\ using $\varpi < \frac 12$ in combination with Markov inequality and 
	\[
	\E\left[\left|\int_{(i-1)\De_n}^{i \De_n} a'_s ds + \int_{(i-1)\De_n}^{i \De_n} \sigma'_s dW_s \right|^q \right] \le C_q \De_n^{\frac q2}
	\]
	for any $q > 0$ by boundedness of $a'$ and $\si'$. 
	
	We finally discuss 
	\begin{align*}
        \P( |{\Delta_i^n Z}| \le K_n , \Delta_i^n R=1) = \P(|\Delta_i^n Z| \le K_n | \Delta_i^n R=1) \P(\Delta_i^n R=1).
   \end{align*}
	Again by properties of the Poisson distribution we have $
	\P(\Delta_i^n R = 1) \le C \Delta_n$ if $n$ is chosen large enough. Hence, we only focus on the former probability which we can bound, using $\xi$ to denote a generic jump size variable, as 
	\begin{align*}
	 &\P\left(\left| \int_{(i-1)\De_n}^{i \De_n} a'_s ds + \int_{(i-1)\De_n}^{i \De_n} \sigma'_s dW_s + \xi \right| \le K_n \right)\\ \le& \P \left(\left| \xi \right| \le 2 K_n \right) + \P\left( \left| \int_{(i-1)\De_n}^{i \De_n} a'_s ds + \int_{(i-1)\De_n}^{i \De_n} \sigma'_s dW_s \right| > K_n \right),
		\end{align*}
	with the second summand being already treated in (\ref{boundexp}). For the first summand we use a Taylor expansion of the jump density $f$ around 0 to get
	\begin{align} \label{ineqf}
	|f(x)| = \left|f(0) + xf^{(1)}(0) + \int_0^x \int_0^y f^{(2)}(z) dz yz  \right| =  \left|\int_0^x \int_0^y f^{(2)}(z) dz yz  \right| \le C x^2
	\end{align}
	for any $x \le 1$, where we have used Condition \ref{condFA}(ii) to obtain $f(0)=f^{(1)}(0)=0$ as well as boundedness of the second derivative. Hence, 
      \begin{align*}
        \P(|{\xi}| \le 2K_n) = \int_{-2K_n}^{2K_n} f(x) dx \le C \int_{-2K_n}^{2K_n} x^2 dx \le C K_n^3 = C \Delta_n^{3\varpi} \le C \Delta_n
      \end{align*}
	for $n$ large enough,	by assumption on $\varpi$. 		
			
  For a fixed \(n\), all of the bounds above hold irrespective of \(i\). Combining these bounds yields (i). 
	
	For part (ii) note first that
	\begin{align*} \P(\De_i^n R = 1) = \De_i^n \La \exp\left(-\De_i^n \La\right)
\end{align*}
from properties of the Poisson distribution. Clearly, 
\begin{align} \label{poissb}
\left| \P(\De_i^n R = 1) - \De_i^n \La \right| \le \De_i^n \La \left(1-\exp\left(-\De_i^n \La\right)\right) \le (\De_i^n \La)^2 \le C \De_n^2
\end{align}
for $n$ large enough, using $1-\exp(-x) \le x$. Then 
	\begin{align*}
	\E[P_n] &= \sum_{i=1}^n \E \left[ \mathbbm{1}_{\{ |\Delta_i^n Z| > K_n \}} \right] = \sum_{i=1}^n \E\left[ \mathbbm{1}_{\{ \Delta_i^n R = 1\}}  \right] + O(n \De_n^2) \\ &= \int_0^{n \De_n} \la_s ds + O(n \De_n^2) = n \De_n \left( \overline \la + o(1) \right) = N_n \left(1 + o(1) \right),
	\end{align*}
	using (i) as well as (\ref{poissb}), $n \De_n \longrightarrow \infty$, and Condition \ref{condFA}(i). We conclude that 
	\[
	\E\left [ \frac{P_n}{N_n} \right] = 1 + o(1). 
	\]
	Then, as in the proof of Lemma \ref{lemPn}(ii), we finally need to show 
	\[
	\Var \left( \frac{P_n}{N_n}  \right) = \frac{\Var(P_n)}{N_n^2} \longrightarrow 0. 
	\]
	By construction, and using properties of the indicator function,
	\begin{align*}
	\Var(P_n) = \sum_{i=1}^n \Var \left( \mathbbm{1}_{\{ |\Delta_i^n Z| > K_n \}} \right) \le \sum_{i=1}^n \E \left[ \mathbbm{1}_{\{ |\Delta_i^n Z| > K_n \}} \right] = N_n \left(1 + o(1) \right),
	\end{align*}
	which finishes the proof as $n \De_n\longrightarrow \infty$. 
	
\subsection{Proof of Theorem \ref{thm2}}
As in the proof of Theorem \ref{thm1} a key step is to identify the main source of the estimation error. We begin with a couple of auxiliary results which we state as a single lemma. 

\begin{lemma}
\label{lemma:e1}
Let \(\xi\) be a random variable with density \(f\) that satisfies Condition \ref{condFA} and set 
\[
\nu_i^n =  \int_{(i-1)\De_n}^{i \De_n} a'_s ds + \int_{(i-1)\De_n}^{i \De_n} \sigma'_s dW_s.
\]
Then the following results hold for $n$ large enough:
\begin{itemize}
	\item[(i)] Let $2 < q < 3$ be arbitrary. Then $\E[|\xi|^{-q}]$ exists and
\begin{equation*}
  \sup_{1 \le i\le n}\left| \E\left[ \left| \nu_i^n + \xi\right|^{-q} \mathbbm{1}_{\{  | \nu_i^n + \xi | > K_n \}} \right] - \E\left[|\xi|^{-q}\right] \right| \le C_q K_n^{3-q}.
\end{equation*}
\item[(ii)] We have
\begin{equation*}
\sup_{1 \le i\le n} \E\left[ \frac{|\nu_i^n|^2}{|\nu_i^n + \xi|^2} \mathbbm{1}_{\{  | \nu_i^n + \xi | > K_n \}} \right] \le C \De_n.
\end{equation*}
\end{itemize}

\end{lemma}

\textit{Step 1:} We prove that the pure jump part is dominating asymptotically, that is 
\begin{equation} \label{decomp1}
    \hat{\beta}_n -\beta = \frac{1}{N_n} \sum_{i=1}^n \sum_{t_{i-1}<s\le t_i}\frac{ \eta_s \Delta Z_s }{\Delta_i^n Z} \mathbbm{1}_{\left\{ |{\Delta_i^n Z} |> K_n \right\}}  + o_p\left(\frac 1{\sqrt{N_n}} \right).
  \end{equation}
Note that Lemma \ref{lemRn} and Slutsky's lemma again allow to replace the scaling $P_n$ with $N_n$, whenever appropriate.  	

For the drift component, with exactly the same proof as in Theorem \ref{thm1}, we know its order to be $\Delta_n^{1-\vpi}$. Then 
\(\sqrt{N_n}\Delta_n^{1-\vpi}=\sqrt{\overline \lambda n \Delta_n^2 \Delta_n^{1-2\vpi}} \longrightarrow 0\) follows easily because of \(n \Delta_n^2\longrightarrow 0\) and \(1-2\vpi> 0\).

For the term involving the continuous martingale component we can show
\[
\frac{1}{N_n} \sum_{i=1}^n  \E\left[ \frac{\int_{(i-1)\De_n}^{i \De_n} \si_s dW_s}{\Delta_i^n Z} \mathbbm{1}_{\{ |\Delta_i^n Z| > K_n \}} \right] = 0
\]
again, as we deal with a sum of martingale differences with mean zero, using the same reasoning as in the proof of Theorem \ref{thm1}. Hence, for any integer $r \ge 2$, 
\begin{align*}
&\E \left[\left| \frac{1}{N_n} \sum_{i=1}^n  \frac{\int_{(i-1)\De_n}^{i \De_n} \si_s dW_s}{\Delta_i^n Z} \mathbbm{1}_{\{ |\Delta_i^n Z| > K_n \}} \right|^2  \right] \\ =& \frac{1}{N_n^2} \sum_{i=1}^n  \E \left[\left| \frac{\int_{(i-1)\De_n}^{i \De_n} \si_s dW_s}{\Delta_i^n Z} \mathbbm{1}_{\{ |\Delta_i^n Z| > K_n \}}   \right|^2  \right] \\ \le& \frac{1}{N_n^2} \sum_{i=1}^n \E \left[\left| \int_{(i-1)\De_n}^{i \De_n} \si_s dW_s \right|^{2r}  \right]^{\frac 1r} \E \left[\left| \Delta_i^n Z\right|^{-\frac{2r}{r-1}} \mathbbm{1}_{\{ |\Delta_i^n Z| > K_n \}}  \right]^{\frac {r-1}r} \\ \le& C_r  \frac{\De_n}{N_n^2} \sum_{i=1}^n \E \left[\left| \Delta_i^n Z\right|^{-\frac{2r}{r-1}} \mathbbm{1}_{\{ |\Delta_i^n Z| > K_n \}} \right]^{\frac {r-1}r}
\end{align*}
where we used both the Hölder inequality and the Burkholder-Davis-Gundy inequality. We will now prove 
\begin{align} \label{ineqmom1}
\E \left[\left| \Delta_i^n Z\right|^{-\frac{2r}{r-1}} \mathbbm{1}_{\{ |\Delta_i^n Z| > K_n \}} \right] \le C_r \left(\De_n +o(1) \right),
\end{align}
uniformly in $i$ and $n$, if we choose $r>3$ and $r$ also large enough such that $\frac{2r}{r-1} \varpi < 1$, which is possible because of $\varpi < \frac 12$. For such a choice of $r$ we can then conclude
\begin{align*}
\E \left[\left| \frac{1}{N_n} \sum_{i=1}^n  \frac{\int_{(i-1)\De_n}^{i \De_n} \si_s dW_s}{\Delta_i^n Z} \mathbbm{1}_{\{ |\Delta_i^n Z| > K_n \}} \right|^2  \right] &\le C_r  \frac{\De_n}{N_n^2} \sum_{i=1}^n \De_n^{\frac {r-1}r} (1+o(1)) \\ &= \frac{n \De_n^2 \De_n^{-\frac 1r}}{N_n^2}(1+o(1)) = o\left(\frac 1{N_n}\right)
\end{align*}
which proves (\ref{decomp1}).

It remains to show (\ref{ineqmom1}) for which we work with a decomposition similarly to the one for Lemma \ref{lemRn}(i): First, we have
\begin{align*}
            &\E\left[\left| \Delta_i^n Z\right|^{-\frac{2r}{r-1}} \mathbbm{1}_{\{ |\Delta_i^n Z| > K_n \}} \mathbbm{1}_{\{ \Delta_i^n R \ge 2 \}} \right] \le K_n^{-\frac{2r}{r-1}}  \P(\Delta_i^n R \ge 2) \le C \Delta_n^{2-\frac{2r}{r-1} \varpi},
          \end{align*}
and then					
          \begin{align*}
            &\E \left[\left| \Delta_i^n Z\right|^{-\frac{2r}{r-1}} \mathbbm{1}_{\{ |\Delta_i^n Z| > K_n \}} \mathbbm{1}_{\{ \Delta_i^n R = 0 \}} \right]  \\
            \le& K_n^{-\frac{2r}{r-1}} \P(|\Delta_i^n Z| > K_n | \Delta_i^n R=0) \P(\Delta_i^n R = 0) \le C_p K_n^{-\frac{2r}{r-1}} \De_n^{p} = C_p \De_n^{p-\frac{2r}{r-1}\varpi}
          \end{align*}
for any $p > 0$, according to (\ref{boundexp}). We then choose $p > 2$ such that this term is asymptotically smaller than the former one.      
Finally, using (\ref{poissb}), $r > 3$, and Lemma \ref{lemma:e1} we obtain				
          \begin{align*}
            &\E\left[ \left| \Delta_i^n Z\right|^{-\frac{2r}{r-1}} \mathbbm{1}_{\{ |\Delta_i^n Z| > K_n \}} \mathbbm{1}_{\{ \Delta_i^n R = 1 \}} \right] \\
            =&  \E\left[\left| \Delta_i^n Z\right|^{-\frac{2r}{r-1}} \mathbbm{1}_{\{ |\Delta_i^n Z| > K_n \}} \middle|  \Delta_i^n R=1 \right] \P(\Delta_i^n R = 1)  \\
            \le& C \Delta_n  \E\left[ \left| \Delta_i^n Z\right|^{-\frac{2r}{r-1}} \mathbbm{1}_{\{ |\Delta_i^n Z| > K_n \}} \middle| \Delta_i^n R=1 \right] \le C_r \De_n \left( \E\left[|\xi|^{{-\frac{2r}{r-1}}}\right] + C_q K_n^{3-\frac{2r}{r-1}} \right). 
          \end{align*}
This proves (\ref{ineqmom1}) due to $\frac{2r}{r-1} \varpi < 1$.

\textit{Step 2:} We show
\begin{equation} \label{decomp2}
    \hat{\beta}_n -\beta = \frac{1}{N_n} \sum_{i=1}^n \sum_{t_{i-1}<s\le t_i}\frac{ \eta_s \Delta Z_s }{\Delta_i^n Z} \mathbbm{1}_{\left\{ |{\Delta_i^n Z} |> K_n \right\}}  \mathbbm{1}_{\left\{ \De_i^n R = 1 \right\}}  + o_p\left(\frac 1{\sqrt{N_n}} \right).
  \end{equation}
	Clearly, we can build upon (\ref{decomp1}). Note first that 	
	\begin{align*}
	\frac{1}{N_n} \sum_{i=1}^n \sum_{t_{i-1}<s\le t_i} \frac{ \eta_s \Delta Z_s }{\Delta_i^n Z} \mathbbm{1}_{\left\{ |{\Delta_i^n Z} |> K_n \right\}} \mathbbm{1}_{\left\{ \De_i^n R =0 \right\}} =0
	\end{align*}
	identically because $\De_i^n R =0$ implies the numerator of the fraction above to vanish. We hence only need to prove 
	\begin{align*}
	\frac{1}{N_n} \sum_{i=1}^n \sum_{t_{i-1}<s\le t_i} \frac{ \eta_s \Delta Z_s }{\Delta_i^n Z} \mathbbm{1}_{\left\{ |{\Delta_i^n Z} |> K_n \right\}} \mathbbm{1}_{\left\{ \De_i^n R \ge 2 \right\}}  = o_p\left(\frac 1{\sqrt{N_n}} \right).
	\end{align*}
	Because both the $\eta_t$ and the jump sizes are independent of all quantities involved, plus $\E[\eta_t] = 0$, we have
	\begin{align*}
	& \E \left[ \left| \frac{1}{N_n} \sum_{i=1}^n \sum_{t_{i-1}<s\le t_i} \frac{ \eta_s \Delta Z_s }{\Delta_i^n Z} \mathbbm{1}_{\left\{ |{\Delta_i^n Z} |> K_n \right\}} \mathbbm{1}_{\left\{ \De_i^n R \ge 2 \right\}}  \right| ^2 \right] \\ =& \frac{1}{N_n^2} \sum_{i=1}^n \E \left[ \left| \sum_{t_{i-1}<s\le t_i} \frac{ \eta_s \Delta Z_s }{\Delta_i^n Z} \mathbbm{1}_{\left\{ |{\Delta_i^n Z} |> K_n \right\}} \mathbbm{1}_{\left\{ \De_i^n R \ge 2 \right\}}  \right| ^2 \right] \\ \le& \frac{1}{N_n^2K_n^2} \sum_{i=1}^n \E \left[ \left| \sum_{t_{i-1}<s\le t_i}  \eta_s \Delta Z_s  \right| ^2 \mathbbm{1}_{\left\{ \De_i^n R \ge 2 \right\}} \right] \\ =& \frac{1}{N_n^2K_n^2} \sum_{i=1}^n \E \left[\eta^2_s \Delta Z^2_s \right] \E \left[ \De_i^n R  \mathbbm{1}_{\left\{ \De_i^n R \ge 2 \right\}} \right] = \frac{\si_\eta^2 \si_\xi^2}{N_n^2K_n^2} \sum_{i=1}^n \E \left[ \De_i^n R  \mathbbm{1}_{\left\{ \De_i^n R \ge 2 \right\}} \right].
	\end{align*}
	Here, $\si_\xi^2$ denotes the second moment of $\xi$. By properties of the Poisson distribution, and using $1-\exp(-x) \le x$, again the right hand side above is bounded by
	\[
	\frac{\si_\eta^2 \si_\xi^2}{N_n^2K_n^2} \sum_{i=1}^n \De_i^n \La \left(1-\exp(-\De_i^n \La) \right) \le C \frac{n \De_n^2}{N_n^2K_n^2} \le C \frac{\De_n^{1-2\varpi}}{N_n}
	\]
for $n$ large enough. This proves (\ref{decomp2}) since $\varpi < \frac 12$. 	

\textit{Step 3:} As a final simplification we show 
\begin{equation} \label{decomp3}
    \hat{\beta}_n -\beta = \frac{1}{N_n} \sum_{i=1}^n \eta_{\tau_i^n}  \mathbbm{1}_{\left\{ \De_i^n R = 1 \right\}}  + o_p\left(\frac 1{\sqrt{N_n}} \right),
  \end{equation}
where $\tau_i^n$ denotes the time of the single jump according to $\De_i^n R = 1$. 
Again, using the previously established (\ref{decomp2}) we have to prove 	
\[
\frac{1}{N_n} \sum_{i=1}^n \left(\sum_{t_{i-1}<s\le t_i}\frac{ \eta_s \Delta Z_s }{\Delta_i^n Z} \mathbbm{1}_{\left\{ |{\Delta_i^n Z} |> K_n \right\}}  - \eta_{\tau_i^n} \right) \mathbbm{1}_{\left\{ \De_i^n R = 1 \right\}} = o_p\left(\frac 1{\sqrt{N_n}} \right),
\]
which we prove in two steps by writing the bracket above as 
\[
\left(\sum_{t_{i-1}<s\le t_i}\frac{ \eta_s \Delta Z_s }{\Delta_i^n Z} - \eta_{\tau_i^n} \right) \mathbbm{1}_{\left\{ |{\Delta_i^n Z} |> K_n \right\}} -\eta_{\tau_i^n}\mathbbm{1}_{\left\{ |{\Delta_i^n Z} | \le  K_n \right\}}. 
\]
For the first term, observe that
\begin{align*}
&\E \left[ \left| \frac{1}{N_n} \sum_{i=1}^n \left(\sum_{t_{i-1}<s\le t_i}\frac{ \eta_s \Delta Z_s }{\Delta_i^n Z}  - \eta_{\tau_i^n} \right) \mathbbm{1}_{\left\{ |{\Delta_i^n Z} |> K_n \right\}} \mathbbm{1}_{\left\{ \De_i^n R = 1 \right\}} \right|^2 \right] \\ =& \frac{\si_\eta^2}{N_n^2} \sum_{i=1}^n  \E \left[ \left|\frac{ J_{\tau_{i}^n}}{\Delta_i^n Z}  - 1 \right|^2 \mathbbm{1}_{\left\{ |{\Delta_i^n Z} |> K_n \right\}}  \mathbbm{1}_{\left\{ \De_i^n R = 1 \right\}} \right] \\ =& \frac{\si_\eta^2}{N_n^2} \sum_{i=1}^n  \E \left[ \left|\frac{ J_{\tau_{i}^n}}{\Delta_i^n Z}  - 1 \right|^2 \mathbbm{1}_{\left\{ |{\Delta_i^n Z} |> K_n \right\}} \middle| \De_i^n R = 1 \right] \P(\De_i^n R=1). 
\end{align*}
For $n$ large enough we have $\P(\De_i^n R = 1) \le C \De_n$ and 
\[
\E \left[ \left|\frac{ J_{\tau_{i}^n}}{\Delta_i^n Z}  - 1 \right|^2 \mathbbm{1}_{\left\{ |{\Delta_i^n Z} |> K_n \right\}} \middle| \De_i^n R = 1 \right] \le C \De_n
\]
according to Lemma \ref{lemma:e1}(ii). Hence,
\begin{align*}
&\E \left[ \left| \frac{1}{N_n} \sum_{i=1}^n \left(\sum_{t_{i-1}<s\le t_i}\frac{ \eta_s \Delta Z_s }{\Delta_i^n Z} \mathbbm{1}_{\left\{ |{\Delta_i^n Z} |> K_n \right\}}  - \eta_{\tau_i^n} \right) \mathbbm{1}_{\left\{ |{\Delta_i^n Z} |> K_n \right\}} \mathbbm{1}_{\left\{ \De_i^n R = 1 \right\}} \right|^2 \right] \\ \le& C \frac{n \De_n^2}{N_n^2} \le C \frac{\De_n}{N_n} = o\left( \frac 1{N_n} \right).
\end{align*}
Finally, 
\begin{align*}
&\E \left[ \left| \frac{1}{N_n} \sum_{i=1}^n \eta_{\tau_i^n} \mathbbm{1}_{\left\{ |{\Delta_i^n Z} | \le  K_n \right\}} \mathbbm{1}_{\left\{ \De_i^n R = 1 \right\}} \right|^2 \right] \\ =& \frac{1}{N_n^2} \sum_{i=1}^n \si^2_\eta  \P \left(|{\Delta_i^n Z} | \le  K_n, \De_i^n R = 1  \right) \le C \frac{n \De_n^{1+3 \varpi}}{N_n^2} \le C \frac{\De_n^{3 \varpi}}{N_n} = o\left( \frac 1{N_n} \right)
\end{align*}
following the proof of Lemma \ref{lemRn}. 

\textit{Step 4:} We finally have to show
\[
 \frac{1}{\sqrt{N_n}} \sum_{i=1}^n \eta_{\tau_i^n}  \mathbbm{1}_{\left\{ \De_i^n R = 1 \right\}} \tol \nN\left(0, \si_\eta^2\right)
\]
which follows directly from the Lyapunov central limit theorem. Setting 
\[
X_i^n = \frac{1}{\sqrt{N_n}} \eta_{\tau_i^n}  \mathbbm{1}_{\left\{ \De_i^n R = 1 \right\}}
\]
it is obvious that $\E[X_i^n] = 0$. We also have 
\begin{align*}
\sum_{i=1}^n \E\left[|X_i^n|^2\right] =& \frac{1}{N_n} \sum_{i=1}^n \E\left[\left|\eta_{\tau_i^n}  \mathbbm{1}_{\left\{ \De_i^n R = 1 \right\}} \right|^2\right] = \frac{1}{N_n} \sum_{i=1}^n\E\left[\eta^2_{\tau_i^n} \middle \vert  \De_i^n R = 1 \right] \P(\De_i^n R = 1) \\ =& \frac{\si_\eta^2 }{N_n} \sum_{i=1}^n \De_i^n \La \exp\left(-\De_i^n \La\right)
\end{align*}
from properties of the Poisson distribution. Using (\ref{poissb}) and Condition \ref{condFA}(i) we obtain
\begin{align*}
\sum_{i=1}^n \E\left[|X_i^n|^2\right] =  \frac{\si_\eta^2 }{N_n} \sum_{i=1}^n \De_i^n \La  + O(\De_n) = \si_\eta^2 \frac 1{n \De_n \overline \la} \int_{0}^{n \De_n} \la_s ds  + O(\De_n)  \longrightarrow \si_\eta^2.
\end{align*}
The same arguments give
\begin{align*}
\sum_{i=1}^n \E\left[|X_i^n|^{2+\de}\right] =& N_n^{-{1+\frac \de 2}} \sum_{i=1}^n\E\left[\left|\eta_{\tau_i^n}\right|^{2+\de} \middle \vert  \De_i^n R = 1 \right] \P(\De_i^n R = 1) \\ \le& C  N_n^{-{1+\frac \de 2}} \sum_{i=1}^n \int_{(i-1)\De_n}^{i \De_n} \la_s ds \le C  N_n^{-{1+\frac \de 2}} \int_{0}^{n \De_n} \la_s ds = o\left( N_n^{-{\frac \de 2}}\right) \longrightarrow 0
\end{align*}
as well. 

\subsubsection{Proof of Lemma \ref{lemma:e1}}
For the proof of (i) note first that
\begin{align} \label{exiq}
\E[|\xi|^{-q}] = \int_{-1}^1 |z|^{-q} f(z) dz +  \int_{\{|z| \ge 1\}}  |z|^{-q} f(z) dz,
\end{align}
and the second integral is clearly finite due to the boundedness of $f$ and $q > 2$. For the first integral, we use $|f(z)| \le C z^2$ on $[-1,1]$ according to (\ref{ineqf}), and then the first integral is finite because of $q < 3$. Hence, $\E[|\xi|^{-q}]$ is finite. 

Now, if we denote the density of $\nu_i^n$ by $g_{i,n}$, and if we write 
\begin{align} \label{defI}
I(\nu) = \E\left[ \left| \nu_i^n + \xi\right|^{-q} \mathbbm{1}_{\{  | \nu_i^n + \xi | > K_n \}} \middle | \nu_i^n = \nu \right]  = \int_{\{|z| > K_n\}} |z|^{-q} f(z - \nu) dz,
\end{align}
then 
\begin{align*}
   &\left| \E\left[ \left| \nu_i^n + \xi\right|^{-q} \mathbbm{1}_{\{  | \nu_i^n + \xi | > K_n \}} \right] - \E\left[|\xi|^{-q}\right] \right|  = \left| \int_{-\infty}^{\infty} I(\nu) g_{i,n}(\nu) d\nu - \E\left[|\xi|^{-q}\right] \right| \\ \le&  \int_{-\infty}^{\infty} \left| I(\nu) - \E\left[|\xi|^{-q}\right] \right|  g_{i,n}(\nu) d\nu.
\end{align*}
To bound the integrand in the previous formula we use a Taylor expansion  around \(z\) and obtain
\begin{align} \label{tayl}
  f(z - \nu) &= f(z) - \nu f^{(1)}(z) + \frac 12 \nu^2 f^{(2)}(z) + T_3(\nu, z)
\end{align}
where we have
\begin{align} \label{trem} \nonumber
 T_3(\nu, z) &= -\int_{z-\nu}^z \int_u^z \int_s^z f^{(3)}(r) dr ds du =   -\int_{z-\nu}^z \int_u^z (r-u) f^{(3)}(r) dr du \\ &= -\frac 12 \int_{z-\nu}^z (r-z+\nu)^2 f^{(3)}(r) dr 
\end{align}
from a twofold application of Fubini's theorem. Hence, the integral in (\ref{defI}) can be written as a sum of four terms, and we discuss each one separately. We first have
\begin{align} \label{aux}
      \left| \int_{\{|z| > K_n\}} |z|^{-q}f(z)  dz - \E\left[|\xi|^{-q}\right] \right|   =  \int_{- K_n}^{K_n} \frac{f(z)}{|z|^q} dz \le  C K_n^{3-q},
 \end{align}
where we again utilize $|f(z)| \le C z^2$ on $[-1,1]$. We also have 
\[
\int_{\{|z| > K_n\}} |z|^{-q} \nu f^{(1)}(z) dz =0, 
\]
for which we use integrability of $f^{(1)}$ as well as the fact that the integrand is an odd function. Furthermore, 
    \begin{align*}
     \int_{\{|z| > K_n\}} \frac 12 \nu^2  |z|^{-q} f^{(2)}(z) dz \le  C \nu^2 \int_{\{|z| > K_n\}} |z|^{-q} dz  =C \nu^2 K_n^{1-q}. 
    \end{align*}
by boundedness of $f^{(2)}$. Finally, 		
\begin{align} \nonumber \label{rembo}
\left| \int_{\{|z| > K_n\}} |z|^{-q} T_3(\nu, z) dz \right| &\le  C \int_{\{|z| > K_n\}} |z|^{-q} \int_{z-\nu}^z (r-z+\nu)^2 \left| f^{(3)}(r) \right| dr dz \\ &\le  C \nu^2 K_n^{-q} \int_{\{|z| > K_n\}} \int_{z-\nu}^z  \left| f^{(3)}(r) \right| dr dz \le C |\nu|^3 K_n^{-q}
\end{align}
using (\ref{trem}), another application of Fubini's theorem and integrability of $f^{(3)}$. To summarize, 
\[
\left| I(\nu) - \E\left[|\xi|^{-q}\right] \right| \le   C \left( K_n^{3-q} + \nu^2 K_n^{1-q} + |\nu|^3 K_n^{-q} \right).
\]
Recalling \(\E[|\nu_i^n|^p] \le C \Delta_n^p\), we finally have
\begin{align*}
   &\left| \E\left[ \left| \nu_i^n + \xi\right|^{-q} \mathbbm{1}_{\{  | \nu_i^n + \xi | > K_n \}} \right] - \E\left[|\xi|^{-q}\right] \right| \le C \left( K_n^{3-q} + \De_n K_n^{1-q} + \De_n^{\frac 32} K_n^{-q} \right) \le C K_n^{3-q}
\end{align*}
for $n$ large enough because of $\varpi < \frac 12$.

The proof of (ii) works in a similar spirit. We have
\begin{align*}
   &\E\left[ \frac{|\nu_i^n|^2}{|\nu_i^n + \xi|^2} \mathbbm{1}_{\{  | \nu_i^n + \xi | > K_n \}} \right] =  \int_{-\infty}^{\infty} L(\nu)  g_{i,n}(\nu) d\nu
\end{align*}
with 
\begin{align*}
L(\nu) =& \E\left[ \frac{|\nu_i^n|^2}{|\nu_i^n + \xi|^2} \mathbbm{1}_{\{  | \nu_i^n + \xi | > K_n \}} \middle | \nu_i^n = \nu \right] = \int_{\{|z| > K_n\}} \frac{\nu^2}{z^2} f(z - \nu) dz \\ =& \int_{\{|z| > K_n\}} \left(\nu^2 S_2(z) + \nu^3 S_3(z) + \nu^4 S_4(z) + S_5(\nu,z) \right) dz, 
\end{align*}
where 
\begin{align*}
S_2(z) = \frac{1}{z^2}f(z), ~ S_3(z) = - \frac{1}{z^2} f^{(1)}(z), ~S_4(z) = \frac{1}{2z^2} f^{(2)}(z), ~S_5(\nu,z) = \frac{\nu^2}{z^2} T_3(\nu,z). 
\end{align*}
Note that the integral over $S_3$ vanishes again due to symmetry. Regarding $S_2$ we first have existence of $\E[\xi^{-2}]$ due to \eqref{exiq}. Then, as in (\ref{aux}) we have
\begin{align*}
      \left| \int_{\{|z| > K_n\}} \frac{f(z)}{z^2} dz - \E\left[|\xi|^{-2}\right] \right|  =   \int_{-K_n}^{K_n} \frac{f(z)}{z^2} dz \le  C K_n. 
    \end{align*}
We can also conclude 	
		\begin{align*}
      \left| \int_{\{|z| > K_n\}} S_4(z) dz \right| \le C \int_{\{|z| > K_n\}} \frac 1{z^2} dz \le C K_n^{-1}.
    \end{align*}
The last component regarding (ii) utilizes the arguments from (\ref{rembo}), as the same reasoning gives 
\begin{align*}
\int_{\{|z| > K_n\}}  \left| S_5(\nu,z) \right| dz \le  C \int_{\{|z| > K_n\}} \frac{\nu^2}{z^2} \int_{z-\nu}^z (r-z+\nu)^2 \left| f^{(3)}(r) \right| dr dz \le C |\nu|^5 K_n^{-2}. 
\end{align*}
In total, 
\[
\left| L(\nu) \right| \le   C \left( \nu^2 (1 + K_n) + \nu^4 K_n^{-1} + |\nu|^5 K_n^{-2} \right),
\]
and from \(\E[|\nu_i^n\|^p] \le C \Delta_n^p\) and $\frac 13 < \varpi < \frac 12$ we finally have
\begin{align*}
   &\E\left[ \frac{|\nu_i^n|^2}{|\nu_i^n + \xi|^2} \mathbbm{1}_{\{  | \nu_i^n + \xi | > K_n \}} \right] \le C \left( \De_n(1+ K_n) + \De_n^2 K_n^{-1} + \De_n^{\frac 32} K_n^{-2} \right) \le C \De_n
\end{align*}
for $n$ large enough. 

\bibliographystyle{chicago}
\bibliography{literatur}

\begin{thebibliography}{}

\bibitem[\protect\citeauthoryear{A\"it-Sahalia, Jacod, and Xiu}{A\"it-Sahalia
  et~al.}{2020}]{aitetal2020a}
A\"it-Sahalia, Y., J.~Jacod, and D.~Xiu (2020).
\newblock Inference on risk premia in continuous-time asset pricing models.
\newblock Working Paper 28140, National Bureau of Economic Research.

\bibitem[\protect\citeauthoryear{A\"it-Sahalia, Kalnina, and Xiu}{A\"it-Sahalia
  et~al.}{2020}]{aitetal2020}
A\"it-Sahalia, Y., I.~Kalnina, and D.~Xiu (2020).
\newblock High-frequency factor models and regressions.
\newblock {\em J. Econometrics\/}~{\em 216\/}(1), 86--105.

\bibitem[\protect\citeauthoryear{Amorino, Jaramillo, and Podolskij}{Amorino
  et~al.}{2024}]{amoetal2024}
Amorino, C., A.~Jaramillo, and M.~Podolskij (2024).
\newblock Optimal estimation of the local time and the occupation time measure
  for an {$\alpha$}-stable {L}\'evy process.
\newblock {\em Mod. Stoch. Theory Appl.\/}~{\em 11\/}(2), 149--168.

\bibitem[\protect\citeauthoryear{Barndorff-Nielsen and
  Shephard}{Barndorff-Nielsen and Shephard}{2004}]{barshe2004}
Barndorff-Nielsen, O.~E. and N.~Shephard (2004).
\newblock Econometric analysis of realized covariation: high frequency based
  covariance, regression, and correlation in financial economics.
\newblock {\em Econometrica\/}~{\em 72\/}(3), 885--925.

\bibitem[\protect\citeauthoryear{Bollerslev, Li, and Todorov}{Bollerslev
  et~al.}{2016}]{boletal2016}
Bollerslev, T., S.~Z. Li, and V.~Todorov (2016).
\newblock Roughing up beta: Continuous versus discontinuous betas and the cross
  section of expected stock returns.
\newblock {\em Journal of Financial Economics\/}~{\em 120\/}(3), 464--490.

\bibitem[\protect\citeauthoryear{Chambers, Mallows, and Stuck}{Chambers
  et~al.}{1976}]{chaetal1976}
Chambers, J.~M., C.~L. Mallows, and B.~W. Stuck (1976).
\newblock A method for simulating stable random variables.
\newblock {\em J. Amer. Statist. Assoc.\/}~{\em 71\/}(354), 340--344.

\bibitem[\protect\citeauthoryear{Cohen}{Cohen}{2013}]{cohen2013}
Cohen, S.~N. (2013).
\newblock A martingale representation theorem for a class of jump processes.

\bibitem[\protect\citeauthoryear{Cont and Tankov}{Cont and
  Tankov}{2004}]{contan2004}
Cont, R. and P.~Tankov (2004).
\newblock {\em Financial modelling with jump processes}.
\newblock Chapman \& Hall/CRC Financial Mathematics Series. Chapman \&
  Hall/CRC, Boca Raton, FL.

\bibitem[\protect\citeauthoryear{Figueroa-L\'opez and Mancini}{Figueroa-L\'opez
  and Mancini}{2019}]{figman2019}
Figueroa-L\'opez, J.~E. and C.~Mancini (2019).
\newblock Optimum thresholding using mean and conditional mean squared error.
\newblock {\em J. Econometrics\/}~{\em 208\/}(1), 179--210.

\bibitem[\protect\citeauthoryear{Jacod and Protter}{Jacod and
  Protter}{2012}]{jacpro12}
Jacod, J. and P.~Protter (2012).
\newblock {\em Discretization of processes}, Volume~67 of {\em Stochastic
  Modelling and Applied Probability}.
\newblock Springer, Heidelberg.

\bibitem[\protect\citeauthoryear{Jacod and Todorov}{Jacod and
  Todorov}{2018}]{jactod2017}
Jacod, J. and V.~Todorov (2018).
\newblock Limit theorems for integrated local empirical characteristic
  exponents from noisy high-frequency data with application to volatility and
  jump activity estimation.
\newblock {\em Ann. Appl. Probab.\/}~{\em 28\/}(1), 511--576.

\bibitem[\protect\citeauthoryear{Li, Todorov, and Tauchen}{Li
  et~al.}{2016}]{lietal2016}
Li, J., V.~Todorov, and G.~Tauchen (2016).
\newblock Inference theory for volatility functional dependencies.
\newblock {\em J. Econometrics\/}~{\em 193\/}(1), 17--34.

\bibitem[\protect\citeauthoryear{Li, Todorov, and Tauchen}{Li
  et~al.}{2017a}]{lietal2017a}
Li, J., V.~Todorov, and G.~Tauchen (2017a).
\newblock Jump regressions.
\newblock {\em Econometrica\/}~{\em 85\/}(1), 173--195.

\bibitem[\protect\citeauthoryear{Li, Todorov, and Tauchen}{Li
  et~al.}{2017b}]{lietal2017b}
Li, J., V.~Todorov, and G.~Tauchen (2017b).
\newblock Robust jump regressions.
\newblock {\em J. Amer. Statist. Assoc.\/}~{\em 112\/}(517), 332--341.

\bibitem[\protect\citeauthoryear{Li, Todorov, Tauchen, and Chen}{Li
  et~al.}{2017}]{lietal2017c}
Li, J., V.~Todorov, G.~Tauchen, and R.~Chen (2017).
\newblock Mixed-scale jump regressions with bootstrap inference.
\newblock {\em J. Econometrics\/}~{\em 201\/}(2), 417--432.

\bibitem[\protect\citeauthoryear{Mancini}{Mancini}{2009}]{mancini2009}
Mancini, C. (2009).
\newblock Non-parametric threshold estimation for models with stochastic
  diffusion coefficient and jumps.
\newblock {\em Scand. J. Stat.\/}~{\em 36\/}(2), 270--296.

\bibitem[\protect\citeauthoryear{Montgomery, Peck, and Vining}{Montgomery
  et~al.}{2001}]{montetal2001}
Montgomery, D.~C., E.~A. Peck, and G.~G. Vining (2001).
\newblock {\em Introduction to linear regression analysis\/} (Third ed.).
\newblock Wiley Series in Probability and Statistics: Texts, References, and
  Pocketbooks Section. Wiley-Interscience, New York.

\bibitem[\protect\citeauthoryear{Mykland and Zhang}{Mykland and
  Zhang}{2006}]{mykzha2006}
Mykland, P.~A. and L.~Zhang (2006).
\newblock A{NOVA} for diffusions and {I}t\^o{} processes.
\newblock {\em Ann. Statist.\/}~{\em 34\/}(4), 1931--1963.

\bibitem[\protect\citeauthoryear{Nolan}{Nolan}{2020}]{nolan2020}
Nolan, J.~P. (2020).
\newblock {\em Univariate stable distributions: models for heavy tailed data}.
\newblock Springer Series in Operations Research and Financial Engineering.
  Springer, Cham.

\bibitem[\protect\citeauthoryear{Schoutens}{Schoutens}{2003}]{schoutens2003}
Schoutens, W. (2003).
\newblock {\em L\'evy Processes in Finance: Pricing Financial Derivatives}.
\newblock Wiley Series in Probability and Statistics. John Wiley \& Sons, Inc.,
  Chichester, England.

\bibitem[\protect\citeauthoryear{Todorov and Bollerslev}{Todorov and
  Bollerslev}{2010}]{todbol2010}
Todorov, V. and T.~Bollerslev (2010).
\newblock Jumps and betas: a new framework for disentangling and estimating
  systematic risks.
\newblock {\em J. Econometrics\/}~{\em 157\/}(2), 220--235.

\end{thebibliography}
\end{document}